\documentclass[a4paper]{article}

\usepackage{amsfonts}
\usepackage{amsmath}
\usepackage{amsthm}\usepackage{amssymb}
\usepackage[mathscr]{eucal}
\usepackage{graphics}
\usepackage[cp866]{inputenc}
\usepackage[english]{babel}
\usepackage{makeidx}
\usepackage{latexsym,amsfonts,amssymb,amsmath,longtable,amsthm}
\input amssym.def
\usepackage{latexsym}
\usepackage{graphicx}
\usepackage[all]{xy}
\usepackage{amsfonts}
\usepackage{amsmath}
\usepackage{mathrsfs}
\usepackage{color}
\usepackage{ulem}
\usepackage[table]{xcolor}

\newtheorem{Theorem}{Theorem}

\newtheorem{Lemma}{Lemma}
\newtheorem{Corollary}{Corollary}
\newtheorem{Remark}{Remark}

 \def\Xint#1{\mathchoice
   {\XXint\displaystyle\textstyle{#1}}%
   {\XXint\textstyle\scriptstyle{#1}}%
   {\XXint\scriptstyle\scriptscriptstyle{#1}}%
   {\XXint\scriptscriptstyle\scriptscriptstyle{#1}}%
   \!\int}
\def\XXint#1#2#3{{\setbox0=\hbox{$#1{#2#3}{\int}$}
     \vcenter{\hbox{$#2#3$}}\kern-.5\wd0}}

\def\dashint{\Xint-}

 \usepackage[pagebackref,
            colorlinks,linkcolor={blue},citecolor={red}]
           {hyperref}

\def\const{{\mathrm{const}}}
\def\pp{{\widehat p}}
\def\0{{\mathbf{0}}}
\def\A{{\mathbf{A}}}
\def\WQ{{\widetilde\Omega}}

\def\RR{{\mathcal{R}}}

\def\Fo{{\mathscr{F}}}

\def\R{{\mathbb R}}
\def\th{{\theta}}
\def\bfeta{\boldsymbol{\theta}}

\def\Ti{{\mathscr{T}}}

\def\H{{\mathscr{H}}}

\def\F{{\mathscr{F}}}

\def\N{{\mathbb{N}}}
\def\Ha{{\mathscr{H}}}

\def\ue{{\mathbf{u}}}
\def\loc{{\mathrm{loc}}}
\def\hh{{\bar h}}

\def\div{\hbox{\rm div}\,}
\def\curl{\hbox{\rm curl}\,}
\def\dist{\hbox{\rm dist}\,}

\def\bmit{\boldsymbol}

 \def\into{\int\limits_}

\def\u{\mathbf u}
\def\w{\mathbf w}

\def\e{\varepsilon}
\def\a{\mathbf a}

\def\n{\mathbf n}

\def\w{\mathbf w}

\def\ve{{\mathbf v}}

\def\div{\hbox{\rm div}\,}
\newcommand{\esssup}{\mathop{\mathrm{ess\,sup}}}
\newcommand{\diam}{\mathop{\mathrm{diam}}}
\newcommand{\meas}{\mathop{\mathrm{meas}}}

\def\F{\mathscr F}

\begin{document}
\title{On the steady Navier--Stokes equations \\in $2D$ exterior domains\footnote{2010 {\it
Mathematical Subject classification}. Primary 76D05, 35Q30;
Secondary 31B10, 76D03; {\it Key words}:   stationary Stokes and
Navier Stokes equations, two--dimensional exterior domains,
boundary value problems.}}

\author{ Mikhail V. Korobkov\footnote{School of Mathematical Sciences,
Fudan University, Shanghai 200433, China; and Voronezh State University,
Universitetskaya pl. 1, Voronezh, 394018, Russia;
korob@math.nsc.ru},  Konstantin Pileckas\footnote{Faculty of
Mathematics and Informatics, Vilnius University, Naugarduko Str.,
24, Vilnius, 03225  Lithuania; konstantinas.pileckas@mif.vu.lt} \,
and Remigio Russo\footnote{Dipartimento di Matematica e Fisica
Universit\`a degli studi della Campania "Luigi Vanvitelli," viale
Lincoln 5, 81100, Caserta, Italy; e-mail:
remigio.russo@unicampania.it}}

 \maketitle

 \date{}

\begin{abstract}
{We study the boundary value problem for the
stationary Navier--Stokes system in two dimensional exterior
domain. We prove that any solution of this problem with finite
Dirichlet integral is uniformly bounded. Also we prove the
existence theorem under zero total flux assumption.}

\end{abstract}

\setcounter{section}{0}

\section{Introduction}

\setcounter{equation}{0}
\setcounter{Theorem}{0}
\setcounter{Lemma}{0}

Let $\Omega$ be an exterior domain in $\R^2$, i.e., 
\begin{equation}\label{Omega}
\Omega={\mathbb R}^2\setminus\bigcup_{i=1}^N\overline\Omega_i,
\end{equation}
where $\Omega_i$ are $N$  pairwise disjoint bounded Lipschitz
domains. The boundary value problem associated with the Navier
Stokes equations in $\Omega$ is to find a solution to  the system
 \begin{equation}
\label{SNS}
\begin{array}{r@{}l}
\nu \Delta{\u}-\u\cdot\nabla\u-\nabla p  & {} ={\bf 0}\qquad \hbox{\rm in } \Omega, \\[2pt]
\div{\u} & {} =0\,\qquad \hbox{\rm in } \Omega ,  \\[2pt]
 \u & {} =\a\qquad\hbox{\rm on }\partial\Omega,
 \end{array}
\end{equation}
with the condition  at infinity
 \begin{equation}
\label{invfrty} \lim_{x\to\infty}\u(x)=\u_0,
\end{equation}
where $\a$ and $\u_0$ are, respectively, an  assigned vector field  on $\partial\Omega$
and a constant vector. Starting from a pioneering
paper by J. Leray \cite{Ler} it is now customary to look for a
solution to  (\ref{SNS}) with   finite Dirichlet integral
 \begin{equation}
\label{SxxS} \into\Omega|\nabla\u|^2dx<+\infty,
\end{equation}
known also as {\it D--solution\/}. As is well known (e.g., \cite{OAL}), such solution is
real--analytic in $\Omega$. Set
 \begin{equation}
\label{SssNS} \mathscr{F}_i=\into{\partial\Omega_i}\a\cdot\n \,ds.
\end{equation}
The existence of a $D$--solution to (\ref{SNS})  has been first
established by J. Leray \cite{Ler} under the assumption
\begin{equation}
\label{SdliS} \mathscr{F}_i=0,\quad i=1\ldots,N.
\end{equation}
To show this,  Leray introduced an elegant argument, known nowadays
as {\it invading domains method\/}, which consists in proving
first that the Navier--Stokes problem
 \begin{equation}
\label{NSDI}
 \begin{array}{rcl}
-\nu \Delta{\bf u}_{k}+\big({\bf u}_{k}\cdot \nabla\big){\bf
u}_{k} +\nabla p_{k} &
= & { 0}\qquad \hbox{\rm in } \Omega_k,\\[4pt]
\div\,{\bf u}_{k} &  = & 0  \qquad \hbox{\rm in } \Omega_k,
\\[4pt]
 {\bf u}_{k} &  = & {\bf a}
 \qquad \hbox{\rm on }\partial\Omega,\\[4pt]
  {\bf u}_{k} &  = & {\bf u}_0
\,\quad \hbox{\rm on }\partial B_k\end{array}
\end{equation}
has a weak solution $ {\bf u}_k$ for every bounded domain
$\Omega_k=\Omega\cap B_k$, $B_k=\{x:|x|<k\}$, $k\gg 1$, and then to
show that the following estimate holds
\begin{equation}
\label{FDIk} \int\limits_{\Omega_k}|\nabla{\bf u}_k|^2dx\le c,
\end{equation}
for some positive constant $c$ independent of $k$. While
(\ref{FDIk}) is sufficient to assure the existence of a subsequence
${\bf u}_{k_l}$ which converges weakly to a  solution~${\bf u}$ of
(\ref{SNS}) satisfying (\ref{SxxS}), it does not give any
information about the behavior at infinity of the velocity ${\bf
u}$\footnote{Indeed, the unbounded function $\log^\alpha|x|$
$(\alpha\in (0,1/2)$) satisfies  (\ref{SxxS}).}, i.e., we do not
know whether ${\bf u}$ satisfies the condition at infinity
\eqref{invfrty}. In 1961 H. Fujita \cite{Fu} recovered, by means
of a different method, Leray's result (see also
\cite[Chapter~XII]{Galdibook}\,). Nevertheless, due to the lack of
a uniqueness theorem, the solutions constructed by Leray and
Fujita are not comparable, even for very small $\nu$.

 Pushing a
little further the argument of Leray \cite{Ler}, A. Russo \cite{AR} showed that the  condition
(\ref{SdliS}) could be extended to the case of "small" (not zero) fluxes by
 \begin{equation}
\label{S234S} \sum_{i=1}^m|\mathscr{F}_i|<2\pi\nu.
\end{equation}

The first existence theorem for (\ref{SNS})--(\ref{invfrty}) is
due to D.R. Smith and R. Finn   \cite{FS}, where it is proved that if
$\u_0\ne{\bf 0}$ and  $|\a-\u_0|$ is sufficiently small, then
there is a $D$--solution to (\ref{SNS}) which converges uniformly
to $\u_0$. This result is particularly meaningful since  it
rules out (at least for small data) for the non--linear Navier--Stokes system
(\ref{SNS})--(\ref{invfrty}) the famous {\it Stokes paradox\/}
which asserts that the equations obtained by  linearization of (\ref{SNS})--(\ref{invfrty})
 \begin{equation}
\label{SNSlin}
\begin{array}{r@{}l}
\nu \Delta{\u} -\nabla p   & {} ={\bf 0}\qquad \hbox{\rm in } \Omega, \\[2pt]
\div{\u} & {} =0\,\qquad \hbox{\rm in } \Omega ,  \\[2pt]
 \u & {} =\a\qquad\hbox{\rm on }\partial\Omega,\\[2pt]
 \displaystyle \lim_{x\to\infty}\u(x)&{} =\u_0,
 \end{array}
\end{equation}
have a solution if and only if
 \begin{equation}
\label{SNPalin}
\into{\partial\Omega}(\a-\u_0)\cdot{\bmit\psi}\,ds={\bf 0},
\end{equation}
for all densities ${\bmit\psi}$ of the simple layer potentials
constant on $\partial\Omega$.  In particular, since
$\int\limits_{\partial\Omega}{\bmit\psi}\ne{\bf 0}$, if $\a$ vanishes  and
$\u_0 $ is a~constant different from zero, then (\ref{SNSlin}) is not solvable.
Moreover, since for the exterior of a ball,  ${\bmit\psi}$ are  the
constant vectors\footnote{More in general, for the exterior of an
ellipsoid  of equation $f(x)=1$,  ${\bmit\psi}={\bmit c}/|\nabla
f|$ for every constant vector ${\bmit c}$ \cite{MRS}.}, a solution
to (\ref{SNSlin})$_{1,2,4}$ satisfies
\begin{equation}
\label{SNaswn} \into0^{2\pi}\u(R,\theta)\,d\theta =2\pi\u_0.
\end{equation}
Of course, by the linearity of the Stokes equations, it is equivalent to
say that a solution to  (\ref{SNSlin})$_{1,2}$ constant on the
boundary and vanishing at infinity does not exist. The situation
is different for the nonlinear problem (\ref{SNS}). The questions whether it admits a solution constant on
$\partial\Omega$ and zero at infinity  is not answered yet, also for
small data.  Nevertheless,   for domains symmetric with respect to
the coordinate axes, i.e.,
$$
(x_1,x_2)\in \Omega\Rightarrow (-x_1,x_2), (x_1,-x_2)\in  \Omega,
$$
in \cite{KRMann} it is showed that a symmetric $D$--solution
\begin{equation}
\label{SymSol}
\begin{array}{ l}
u_1(x_1,x_2)=-u_1(-x_1,x_2)=u_1(x_1,-x_2) \\[2pt]
u_2(x_1,x_2)=u_2(-x_1,x_2)=-u_2(x_1,-x_2),
 \end{array}
\end{equation}
to (\ref{SNS}), uniformly vanishing at infinity, exists under the
only natural assumption that $\a$ satisfies (\ref{SymSol}) and natural regularity conditions. Note
that (\ref{SymSol}) meets the mean property (\ref{SNaswn}) with
$\u_0={\bf 0}$.

The problem  of the asymptotic behavior at infinity of an arbitrary
$D$--solution $(\u,p)$ to (\ref{SNS})$_{1,2}$ was
tackled by D.~Gilbarg \& H. Weinberger \cite{GW1}--\cite{GW2} and
C.~Amick \cite{Amick}. In \cite{GW2} it is shown that
\begin{equation}\label{prc}
p-p_0=o(1)\qquad\mbox{ as }r\to\infty,
\end{equation}
i.e., pressure has a limit at infinity (one can choose, say, $p\to
0$\,), and
\begin{equation}
\label{SzzS}
\begin{array}{ l}
 \u(x) =o(\log^{1/2} r),\\[2pt]
 \omega  =o(r^{-3/4}\log^{1/8}r),\\[2pt]
\nabla\u(x)=o(r^{-3/4}\log^{9/8}r),
 \end{array}
\end{equation}
where
$$
\omega=\partial_1u_2-\partial_2u_1
$$
is the corresponding  vorticity. If, in addition,  $\u$ is
bounded, then  there is a constant vector $\u_\infty$ such that
 \begin{equation}
\label{GGWW}
\displaystyle\lim_{r\to+\infty}\into0^{2\pi}|\u(r,\theta)-\u_\infty|^2d\theta=0,\end{equation}
and
\begin{equation}
\label{Sz123zS}
\begin{array}{l}
 \omega  =o(r^{-3/4}),\\[2pt]
\nabla\u(x)=o(r^{-3/4}\log r).
 \end{array}
\end{equation}
Here if $\u_\infty={\bf 0}$, then $\u=o(1)$.  Moreover, in
\cite{ARPJ} it is proved  that 
\begin{equation}
\label{SzzasS} \nabla p  =O(r^{\epsilon-1/2})\end{equation} 
for every positive $\epsilon$.

In
\cite{Amick} it is proved  that  if $\u$ vanishes on the boundary,
then $\u$ is bounded and, as a consequence,  satisfies
(\ref{GGWW}), (\ref{Sz123zS}). However, in this last case the
solution could  tend  to zero at infinity and even be the
trivial one. This possibility was excluded  by  Amick
\cite{Amick} (Section 4.2) for the solution obtained by the Leray
method, for symmetric  with respect to the $x_2$--axis
(say) domains, i.e, $(x_1,x_2)\in \Omega\Rightarrow (x_1,-x_2)\in
\Omega$.  This result is remarkable as the first step
to exclude  the non--linear Stokes paradox for every $\nu$, at
least for axisymmetric domains. For such  kind of domains
the existence of a $D$--solution to (\ref{SNS}) is established   in
\cite{kprplane} only under the symmetry hypothesis
$a_1(x_1,x_2)=a_1(x_1,-x_2)$, $a_2(x_1,x_2)=-a_2(x_1,-x_2)$.

Despite the efforts of many researchers (see, $e.g$, the reference
in \cite{Galdibook}) several relevant  problems remain open, among
which: existence of a   solution    to (\ref{SNS}) for arbitrary
fluxes $\mathscr{F}_i$, its uniqueness (for small data);  the boundedness of
a~$D$--solutions (in the case of non-homogeneous boundary conditions), its uniform convergence to  $\u_\infty\ne{\bf
0}$\footnote{ By a remarkable result of L.I. Sazonov
\cite{Sazonov}, this ensures  that the solution behaves at
infinity as that of the linear Oseen equations (see also
\cite{GaldiH} and \cite{Galdibook})} and the relation  between
$\u_\infty$ and  $\u_0$; more precise asymptotic behavior of
$\nabla p$ and the derivatives of $\u$.

The present paper is devoted  to some of the above issues. The
first main result is as follows.

\begin{Theorem}
\label{T1} {\sl Let $\u$ be a solution to the Navier--Stokes
system
 \begin{equation}
\label{SNS-1} \left\{\begin{array}{r@{}l}
-\nu \Delta{\u}+\u\cdot\nabla\u+\nabla p  & {} ={\bf 0}\qquad \hbox{\rm in } \Omega, \\[2pt]
\div{\u} & {} =0\,\qquad \hbox{\rm in } \Omega
 \end{array}\right.
\end{equation}
in the exterior domain~$\Omega\subset\R^2$.
Suppose
\begin{equation}
\label{dom4} \int\limits_{\Omega}|\nabla\u|^2dx<\infty.
\end{equation}
Then $\u$ is uniformly bounded in $\Omega_0=\R^2\setminus B_{R_0}$, i.e.,
\begin{equation}
\label{dom2} \sup\limits_{x\in\Omega_0}|\u(x)|<\infty,
\end{equation}
where $B_{R_0}$ is a~disk with
sufficiently large radius: $\frac12 B_{R_0}\Supset\partial\Omega$.}
 \end{Theorem}

Using the above--mentioned results of  D.~Gilbarg and
H.~Weinberger, we obtain immediately

\begin{Corollary}\label{GW-cor}
{\sl Let $\u$ be a $D$-solution to the Navier--Stokes
system~(\ref{SNS-1}) in a neighbourhood of infinity.  Then the
asymptotic properties~(\ref{prc}), (\ref{GGWW})--(\ref{Sz123zS})
hold.}
\end{Corollary}

Using the results of the above--mentioned paper of
Amick~\cite{Amick}, we could say something more about asymptotic
properties of~$D$-solutions in the case of zero total flux, i.e.,
when

\begin{equation}
\label{dom5} \int\limits_{\partial\Omega}\a\cdot \n\,ds=0,
\end{equation}

\begin{Corollary}\label{Am-cor}
{\sl Let $\u$ be a $D$-solution to the Navier--Stokes
problem~(\ref{SNS-1}) in an~exterior domain~$\Omega\subset\R^2$ with
zero total flux condition~(\ref{dom5}). Then in addition to the
properties of Theorem~\ref{T1} and Corollary~\ref{GW-cor}, the
total head pressure $\Phi=p+\frac12|\u|^2$ and the absolute value
of the velocity $|\u|$ have the uniform limit at infinity, i.e.,
\begin{equation}
\label{dom6} |\u(r,\th)|\to|\u_\infty|\qquad\mbox{ as }\
r\to\infty,
\end{equation}
where $\u_\infty$ is a constant vector from the
condition~(\ref{GGWW}).}
\end{Corollary}

Let us note that formally Amick \cite{Amick} established~(\ref{dom6})  under the~stronger assumption
\begin{equation}
\label{dom5-} \a\equiv 0.
\end{equation}
But really his argument for~(\ref{dom6}) cover the more
general case~(\ref{dom5}) as well. Indeed, the main tool
in~\cite{Amick} was the use of the~auxiliary function
$\gamma=\Phi-\omega\psi$, where $\psi$ is a stream function:
$\nabla\psi=\u^\bot=(u_2,-u_1)$. This auxiliary function~$\gamma$
has remarkable monotonicity properties: it is monotone along
 level sets of the~vorticity~$\omega=c$ and  vice versa~--
the~vorticity is monotone along level sets~$\gamma=c$. But, of course,
the stream function~$\psi$ (and, consequently, the corresponding
auxiliary function~$\gamma$\,) could be well defined in the
neighbourhood of infinity under the more general case~(\ref{dom5})
instead of~(\ref{dom5-}). Furthermore,  Amick also proved
that under the~conditions of Corollary~\ref{Am-cor}, the convergence 
\begin{equation}
\label{dom8} \gamma(r,\theta)\to\frac12|\u_\infty|\qquad\mbox{ as
}\ r\to\infty
\end{equation}
holds uniformly with respect to $\theta$.

The second result of the paper concerns the existence of solutions to the non-homogeneous boundary value problem~\eqref{SNS}.

\begin{Theorem}
\label{T2} {\sl Let $\Omega\subset\R^2$ be an exterior domain with
$C^2$-smooth boundary. Suppose that $\a\in
W^{1/2,2}(\partial\Omega)$ and the equality~(\ref{dom5}) holds,
i.e., the total flux is zero. Then there exists a~$D$-solution
$\u$ to the Navier--Stokes boundary value problem~$(\ref{SNS})$. }
\end{Theorem}

This theorem shows also that the asymptotic results of Corollaries~\ref{GW-cor},~\ref{Am-cor} \ and \
(\ref{dom8}) have meaning and are not just a figment of the
imagination.  

Note, that the existence theorem
 for the steady Navier--Stokes problem in {\it three dimensional}
 exterior axially symmetric domains (with axially symmetric data) was proved in the recent
 paper~\cite{KPRMA} without any conditions on fluxes~$\mathscr{F}_i$.

\section{Notations and preliminaries}
\setcounter{equation}{0}

By {\it a domain} we mean an open connected set. We use standard
notations for function spaces: $W^{k,q}(\Omega)$,
$W^{\alpha,q}(\partial\Omega)$, where $\alpha\in(0,1),
k\in{\mathbb N}_0, q\in[1,+\infty]$. In our notation we do not
distinguish function spaces for scalar and vector valued
functions; it is clear from the context whether we use scalar or
vector (or tensor) valued function spaces.

For $q\ge1$ denote by $D^{k,q}(\Omega)$ the set of functions $f\in
W^{k,q}_{\loc}(\Omega)$ such that
$\|f\|_{D^{k,q}(\Omega)}=\|\nabla^k f\|_{L^q(\Omega)}<\infty$.
Further, $D^{1,2}_0(\Omega)$ is the closure of the set of all
smooth functions having compact supports in $\Omega$ with respect
to the norm $\|\,\cdot\,\|_{D^{1,2}(\Omega)}$, \, and
$H(\Omega)=\{{\bf v}\in D^{1,2}_0(\Omega):\, \div {\bf v}=0\}$; \
$D^{1,2}_\sigma(\Omega):=\{{\bf v}\in D^{1,2}(\Omega):\, \div {\bf
v}=0\}$.

\section{Boundedness of general~$D$-solutions: proof of Theorem~\ref{T1}.}
\setcounter{equation}{0}

Suppose the assumptions of Theorem~\ref{T1} are fulfilled. By
classical regularity results for $D$-solutions to Navier--Stokes
system,  the function~$\u$ is uniformly bounded on each bounded
subset of the set~$\Omega_0=\R^2\setminus B_{R_0}$; moreover, $\u$ is analytical in~$\Omega_0$.
By results of~\cite{GW2},  pressure  is uniformly bounded in~$\Omega_0$:
\begin{equation}
\label{dom10} \sup\limits_{x\in\Omega_0}|p(x)|\le C<+\infty.
\end{equation}
Suppose that the assertion~(\ref{dom2}) of the~Theorem
is false. Then there exists a sequence of points $x_k\in\Omega_0$
such that
\begin{equation}
\label{dom11} |x_k|\to+\infty\qquad\mbox{ \ and \quad
}|\u(x_k)|\to+\infty.
\end{equation} This means, by virtue of~(\ref{dom10}), that
\begin{equation}
\label{dom11'} \Phi(x_k)\to+\infty,
\end{equation} where
$\Phi=p+\frac12|\u|^2$ is the total head pressure.

Since $\u$ is a $D$-solution,
$\int\limits_{\Omega_0}|\nabla\u|^2dx<\infty$, by standard arguments
there exists an increasing sequence on numbers $R_m<R_{m+1}$ such
that $R_m\to\infty$ and
\begin{equation}
\label{dom12} \int\limits_{{C}_{R_m}}|\nabla \u|\,ds\to 0,
\end{equation}
where ${C}_R:=\{x\in\R^2:|x|=R\}$. It implies that
\begin{equation}
\label{dom13} \sup\limits_{x\in {C}_{R_m}}|\u(x)-\bar\u_m|\to 0,
\end{equation}
here $\bar\u_m$ is the mean value of $\u$ on the circle~${C}_{R_m}$. Indeed, for any component $u_j$ of $\u$,  by mean value theorem,  there exists a point $\theta_j^*\in [0,2\pi)$ such that
$$
u_j(R_m, \theta_{j}^{*})=(2\pi)^{-1}\int\limits_{0}^{2\pi}u_j(R_m, \theta)d\theta=\bar u_{jm},\quad j=1,2,
$$
and
$$
|u_j(R_m, \theta)-\bar u_{jm}|=|u_j(R_m, \theta)-u_j(R_m, \theta_j^*)|\leq \int\limits_{\theta_j^*}^{\theta}\big|\frac{\partial u_j}{\partial\theta}\big|d\theta\leq \int\limits_{C_{R_m}}|\nabla \u|\,ds\rightarrow 0.
$$

Since $\Phi$ satisfies the maximum principle (see,
e.g.,\cite{GW2}\,), in particular, for any subdomain
$\Omega_{m_1,m_2}=\{x:R_{m_1}<|x|<R_{m_2}\},$ with $\partial\Omega_{m_1,m_2}=C_{R_{m_1}}\cup
{C}_{R_{m_2}}$ we have
$$\sup\limits_{x\in\Omega_{m_1,m_2}}\Phi(x)=\sup\limits_{x\in {C}_{R_{m_1}}\cup {C}_{R_{m_2}}}\Phi(x).
$$
Relations  (\ref{dom11}), (\ref{dom13}) imply that
$|\bar\u_m|\to+\infty$; consequently, by \eqref{dom10}, \eqref{dom11'}, \eqref{dom13},
$$\inf\limits_{x\in {C}_{R_{m}}}\Phi(x)\to+\infty.
$$
Then we could assume without loss of generality (choosing a subsequence) that
\begin{equation}
\label{dom14} \sup\limits_{x\in
{C}_{R_{m}}}\Phi(x)<\inf\limits_{x\in {C}_{R_{m+1}}}\Phi(x).
\end{equation}
Recall that by the classical Morse--Sard Theorem (see, e.g.,
\cite{Hir}\,), applied to the analytical function~$\Phi$, for
almost all values $t\in\Phi(\Omega_0)$ the level set $\{\Phi=t\}$
contains no critical points, i.e., $\nabla \Phi(x)\ne0$ if
$x\in\Omega_0$ and $\Phi(x)=t$. Further such values are
called~{\it regular}. Take arbitrary regular value
$t>t_*=\sup\limits_{x\in {C}_{R_{1}}\cup {C}_{R_0}}\Phi(x)$. Then
by the~implicit function theorem the level set
$\{x\in\Omega_0:\Phi(x)=t\}$ consists of a family of disjoint
smooth curves which are separated (by construction) both from
infinity and from the boundary~$\partial\Omega_0={C}_{R_0}$. Of
course, this implies that every connected component of this level
set $\{\Phi=t\}$ is homeomorphic to a~circle. Let us call these
components {\it quasicircles}. By obvious geometrical arguments,
for every regular $t>t_*$ there exists at least one quasicircle
$S$ separating $C_{R_1}$ from infinity, i.e., $C_{R_1}$ is
contained in the bounded connected component of the open
set~$\R^2\setminus S$. Because of the~maximum principle, such
quasicircle is unique, and we will denote it by~$S_t$.

For $t_*<\tau<t$ let $\Omega_{\tau,t}$ be a domain with
$\partial\Omega_{\tau,t}=S_\tau\cup S_t$. Integrating the identity
\begin{equation}
\label{dom15} \Delta\Phi=\omega^2+\frac1\nu\div(\Phi\u)
\end{equation}
over  $\Omega_{\tau,t}$, we obtain
\begin{equation}
\label{dom16}
\begin{array}{lcr}
\int\limits_{S_t}|\nabla\Phi|\,ds-\int\limits_{S_\tau}|\nabla\Phi|\,ds=\int\limits_{\Omega_{\tau,t}}\omega^2dx+\frac1\nu\int\limits_{S_t}\Phi \u\cdot\n\,ds-\frac1\nu\int\limits_{S_\tau}\Phi \u\cdot\n\,ds\\
\\
=\int\limits_{\Omega_{\tau,t}}\omega^2dx+\frac1\nu(t-\tau)\Fo,
\end{array}
\end{equation}
where $\Fo=\int\limits_{{C}_{R_0}}\u\cdot\n$ is the total flux. Notice that by construction the unit normal $\n$ to the level set $S_t=\{x: \Phi(x)=t\}$ is equal to $\frac{\nabla \Phi}{|\nabla\Phi|}$,
so that $\nabla\Phi\cdot\n=|\nabla\Phi|$ on $S_t$; \ analogously, $\nabla\Phi\cdot\n=-|\nabla\Phi|$ on $S_\tau$.
The further proof splits into two cases.

{\sc Case~I}. The total flux in not zero: $\Fo\ne 0$. First suppose that $\Fo>0$.
Then from~(\ref{dom16}) \ (fixing $\tau$ and taking a~big~$t$\,)  we
obtain
\begin{equation}
\label{dom17} C_1t\le\int\limits_{S_t}|\nabla\Phi|\,ds\le C_2 t
\end{equation}
for sufficiently large~$t$ and for some positive
constants~$C_1,C_2$ (not depending on~$t$). Denote by $\RR$ the
set of all regular values $t>t_*$, and put
$$E_t:=\bigcup\limits_{\tau\in[t,2t]\cap\RR}S_\tau.$$
Applying the classical Coarea formula
$$
\int\limits_{E_t}f\,|\nabla \Phi|\,dx=\int\limits_t^{2t}\biggl(\int\limits_{S_\tau}f\,ds\biggr)\,d\tau
$$
for $f=|\omega|$ and for $f=|\nabla\Phi|$ we obtain
\begin{equation}
\label{dom18}
\begin{array}{lcr}
\int\limits_t^{2t}\biggl(\int\limits_{S_\tau}|\omega|\,ds\biggr)\,d\tau=\int\limits_{E_t}|\omega|\cdot|\nabla
\Phi|dx\le\biggl(\int\limits_{E_t}|\nabla\Phi|^2dx\biggr)^{\frac12}
\biggl(\int\limits_{E_t}\omega^2dx\biggr)^{\frac12}\\
\\
=\biggl(\int\limits_t^{2t}\biggl(\int\limits_{S_\tau}|\nabla\Phi|\,ds\biggr)\,d\tau\biggr)^{\frac12}
\biggl(\int\limits_{E_t}\omega^2dx\biggr)^{\frac12} \le\e t,
\end{array}
\end{equation}
where $\e\to0$ as $t\to\infty$ (we used here~(\ref{dom17}) and the
assumption that  the Dirichlet integral is finite).
From
(\ref{dom18}) and from the mean value theorem it
follows that there exists a value $\tau\in[t,2t]\cap\RR$ such that
\begin{equation}
\label{dom19} \int\limits_{S_\tau}|\omega|ds\le 2\e.
\end{equation}
Since the pressure is uniformly bounded (see \eqref{dom10}),  we conclude that
$|\u|\sim\sqrt{2\tau}$ on $S_\tau$ for large~$\tau$, therefore,
using the identity
$$\nabla\Phi=-\nu\nabla^\bot\omega+\omega\u^\bot,$$
we obtain
\begin{equation}
\label{dom20}
\int\limits_{S_\tau}|\nabla\Phi|\,ds=\int\limits_{S_\tau}\omega\u^{\bot}\cdot\n\,ds\le2\sqrt{\tau}\,\int\limits_{S_\tau}|\omega|\,ds\le
4\sqrt{\tau}\e
\end{equation}
(the integral of $\nabla^\bot\omega\cdot\n=\curl \omega\cdot \n$ over the closed curve $S_\tau$ is equal to zero).

The last estimate contradicts the first inequality in~(\ref{dom17}). Thus, if $\Fo>0$, then the~assumption \eqref{dom11} is false and the solution $\u$ is uniformly bounded.

Let $\Fo<0$. Writing relation \eqref{dom16} in the form
\begin{equation}
\label{dom16'}
\int\limits_{S_t}|\nabla\Phi|\,ds=\int\limits_{S_\tau}|\nabla\Phi|\,ds+\int\limits_{\Omega_{\tau,t}}\omega^2dx+\frac1\nu(t-\tau)\Fo,
\end{equation}
we immediately see that for large $t$ the right-hand side becomes negative, while the left-hand side is positive for all $t$. We again obtain a contradiction to assumption \eqref{dom11}.
Thus, the proof for the case $\Fo\neq 0$ is complete.

\

{\sc Case~II}. The total flux is zero: $\Fo= 0$.
Then formula~(\ref{dom16}) takes the form
\begin{equation}
\label{dom16-}
\int\limits_{S_t}|\nabla\Phi|\,ds=\int\limits_{S_\tau}|\nabla\Phi|\,ds+\int\limits_{\Omega_{\tau,t}}\omega^2dx.
\end{equation}
From the last identity it follows that
$\int\limits_{S_t}|\nabla\Phi|\,ds$ is a bounded increasing function,
i.e., it has a~finite positive limit, in particular,
\begin{equation}
\label{dom17-} C_1\le\int\limits_{S_t}|\nabla\Phi|\,ds\le C_2
\end{equation}
for sufficiently large~$t$ and for some positive
constants~$C_1,C_2$ (independent of~$t$). Applying the Coarea formula,  we obtain now
\begin{equation}
\label{dom18-}
\begin{array}{lcr}
\int\limits_t^{2t}\biggl(\int\limits_{S_\tau}|\omega|\,ds\biggr)\,d\tau=\int\limits_{E_t}|\omega|\cdot|\nabla
\Phi|dx\le\biggl(\int\limits_{E_t}|\nabla\Phi|^2dx\biggr)^{\frac12}\cdot
\biggl(\int\limits_{E_t}\omega^2dx\biggr)^{\frac12}\\
\\
=\biggl(\int\limits_t^{2t}\biggl(\int\limits_{S_\tau}|\nabla\Phi|\,ds\biggr)\,d\tau\biggr)^{\frac12}\cdot
\biggl(\int\limits_{E_t}\omega^2dx\biggr)^{\frac12} \le\e \sqrt{t},
\end{array}
\end{equation}
where $\e\to0$ as $t\to\infty$. From
(\ref{dom18-}) and from the mean value theorem the existence of a~value $\tau\in[t,2t]\cap\RR$
follows such that
\begin{equation}
\label{dom19-} \int\limits_{S_\tau}|\omega|\,ds\le \e\frac2{\sqrt{\tau}}.
\end{equation}
As in the Case I we have
$|\u|\sim\sqrt{2\tau}$ on $S_\tau$. Therefore, integrating again the~identity
$$\nabla\Phi=-\nu\nabla^\bot\omega+\omega\u^\bot,$$
we obtain
\begin{equation}
\label{dom20-}
\int\limits_{S_\tau}|\nabla\Phi|\,ds=\int\limits_{S_\tau}\omega\u^{\bot}\cdot\n\,ds\le2\sqrt{\tau}\,\int\limits_{S_\tau}|\omega|dx\le
4\e.
\end{equation}
The last estimate is in
contradiction with the first inequality in~(\ref{dom17-}). Therefore, in the case  $\Fo=0$  assumption \eqref{dom11} is again false and the solution $\u$ is uniformly bounded.
Theorem~\ref{T1} is proved.

\

\section{The existence theorem: proof of Theorem~\ref{T2}.}
\setcounter{equation}{0}

Here we need some preliminary results on real analysis and
topology.

\subsection{On Morse-Sard and Luzin N-properties of Sobolev
functions from $W^{2,1}$} \label{SMS}

Let us recall some classical differentiability properties of
Sobolev functions.

\begin{Lemma}[see Proposition~1 in \cite{Dor}]
\label{kmpThDor} {\sl Let  $\psi\in W^{2,1}(\R^2)$. Then the
function~$\psi$ is continuous and there exists a set $A_{\psi}$
such that $\mathfrak{H}^1(A_{\psi})=0$, and the function $\psi$ is
differentiable (in the classical sense) at each $x\in\R^2\setminus
A_{\psi}$. Furthermore, the classical derivative at such points
$x$ coincides with $\nabla\psi(x)=\lim\limits_{r\to 0}
\dashint_{B_r(x)}{ \nabla\psi}(z)dz$, and \ $\lim\limits_{r\to
0}\dashint\nolimits_{B_r(x)}|\nabla
\psi(z)-\nabla\psi(x)|^2dz=0$.}
\end{Lemma}

Here and henceforth we denote by $\mathfrak{H}^1$ the
one-dimensional Hausdorff measure, i.e.,
$\mathfrak{H}^1(F)=\lim\limits_{t\to 0+}\mathfrak{H}^1_t(F)$,
where $\mathfrak{H}^1_t(F)=\inf\{\sum\limits_{i=1}^\infty {\rm
diam} F_i:\, {\rm diam} F_i\leq t, F\subset
\bigcup\limits_{i=1}^\infty F_i\}$.

The next theorem have been proved recently by J. Bourgain,
M.~Korobkov and J. Kristensen \cite{korob} (see also \cite{BKK2}
for a multidimensional case).

\begin{Theorem}
\label{kmpTh1.1}{\sl  Let  ${\mathcal D}\subset\R^2$ be a bounded
domain with Lipschitz boundary and  $\psi\in W^{2,1}({\mathcal
D})$. Then

{\rm (i)} $\mathfrak{H}^1(\{\psi(x)\,:\,x\in\bar{\mathcal
D}\setminus A_\psi\,\,\&\,\,\nabla \psi(x)=0\})=0$;

{\rm (ii)} for every $\varepsilon>0$ there exists $\delta>0$ such
that for any set $U\subset \bar{\mathcal D}$ with
$\mathfrak{H}^1_\infty(U)<\delta$ the inequality
$\mathfrak{H}^1(\psi(U))<\varepsilon$ holds;

{\rm (iii)} for $\mathfrak{H}^1$--almost all $y\in
\psi(\bar{\mathcal D})\subset \mathbb{R}$ the preimage
$\psi^{-1}(y)$ is a finite disjoint family of $C^1$--curves $S_j$,
$ j=1, 2, \ldots, N(y)$. Each $S_j$ is either a cycle in
${\mathcal D}$ $($i.e., $S_j\subset{\mathcal D}$ is homeomorphic
to the unit circle $\mathbb{S}^1)$ or it is a simple arc with
endpoints on $\partial{\mathcal D}$ $($in this case $S_j$ is
transversal to $\partial{\mathcal D}\,)$. }
\end{Theorem}

\subsection{Some facts from topology}
\label{Kronrod-s}

We shall need  some topological definitions and results. By  {\it
continuum} we mean a compact connected set. We understand
connectedness  in the sense of general topology. A~subset of a
topological space is called {\it an arc} if it is homeomorphic to
the unit interval~$[0,1]$.

Let us shortly present  some results from the classical paper of
A.S.~Kron-\\rod~ \cite{Kronrod} concerning level sets of
continuous functions. Let ${Q}=[0,1]\times[0,1]$ be a square in
$\mathbb{R}^2$ and let $f$ be a continuous function  on ${Q}$.
Denote by $E_t$ a level set of the function $f$, i.e.,
$E_t=\{x\in{Q}: f(x)=t\}$. A component $K$  of the level set $E_t$
containing a point $x_0$ is a maximal connected subset of $E_t$
containing $x_0$. By $T_f$ denote a family of all connected
components of level sets of~$f$. It was established
in~\cite{Kronrod} that $T_f$ equipped by a natural
topology\footnote{The convergence in $T_f$ is defined as follows:
$T_f\ni C_i\to C$ iff $\sup\limits_{x\in C_i}\dist(x,C)\to0$.} is
a one-dimensional topological tree\footnote{A locally connected
continuum $T$ is called~{\it a topological tree}, if it does not
contain a curve homeomorphic to a circle, or, equivalently, if any
two different points of~$T$ can be joined by a unique arc. This
definition implies that $T$ has topological dimension~1.}.
Endpoints of this tree\footnote{A point of a continuum~$K$ is
called an {\it endpoint of $K$} (resp., {\it a branching point
of~$K$}) if its topological index equals~1 (more or equal to~$3$
resp.). For a~topological tree~$T$ this definition is equivalent
to the~following: a point $C\in T$ is an~endpoint of~$T$ (resp., {
a branching point of~$T$}), if the set $T\setminus\{C\}$ is
connected (resp., if $T\setminus\{C\}$ has more than two connected
components).} are the components~$C\in T_f$ which do not
separate~$Q$, i.e., $Q\setminus C$ is a connected set. Branching
points of the tree are the components $C\in T_f$ such that
$Q\setminus C$ has more than two connected components (see
\cite[Theorem 5]{Kronrod}). By results of \cite[Lemma~1]{Kronrod},
the set of all branching points of~$T_f$ is at most countable. The
main property of a tree is that any two points could be joined by
a unique arc. Therefore, the same is true for~$T_f$.

\begin{Lemma} [see Lemma~13 in \cite{Kronrod}]
\label{kmpLem6} If $f\in C(Q)$, then for any two different points
$A\in T_f$ and $B\in T_f$, there exists a unique arc
$J=J(A,B)\subset T_f$ joining $A$ to $B$. Moreover, for every
inner point $C$ of this arc the points $A,B$ lie in  different
connected components of the set $T_f\setminus\{C\}$.
\end{Lemma}

We can reformulate the above Lemma in the following equivalent
form.

\begin{Lemma} \label{kmpLem7}{\sl If  $f\in C(Q)$, then for
any two different points $A,B\in T_f$, there exists a~continuous
injective function $\varphi:[0,1]\to T_f$ with the properties

{\rm (i)} $\varphi(0)=A$, $\varphi(1)= B$;

{\rm (ii)} for any $t_0\in[0,1]$,
$$
\lim\limits_{[0,1]\ni t\to t_0}\sup\limits_{x\in
\varphi(t)}\dist(x,\varphi(t_0))\to0;
$$

{\rm (iii)}  for any $t\in(0,1)$ the sets $A,B$ lie in different
connected components of the set \ $Q\setminus\varphi(t)$.}
\end{Lemma}

\begin{Remark}
\label{kmpRem2} {\rm If in  Lemma~\ref{kmpLem7} $f\in W^{2,1}(Q)$,
then by Theorem~\ref{kmpTh1.1}~(iii), there exists a dense subset
$E$ of $(0,1)$ such that $\varphi(t)$ is a $C^1$-- curve for every
$t\in E$. Moreover, $\varphi(t)$ is either a cycle   or a simple
arc with endpoints on $\partial Q$.}
\end{Remark}

\begin{Remark}
\label{kmpRem1.2} {\rm All results of
Lemmas~\ref{kmpLem6}--\ref{kmpLem7} remain valid for level sets of
continuous functions $f:\overline\Omega_0\to\R$, where
$\overline\Omega_0\subset\R^2$ is a compact set homeomorphic to
the unit square $Q=[0,1]^2$.}
\end{Remark}
\

\subsection{Leray's argument ``reductio ad absurdum'' }
\label{poet}

Consider the Navier--Stokes problem~(\ref{SNS})  in the
$C^2$-smooth
 exterior domain $\Omega\subset\R^2$ defined
by~(\ref{Omega}).  Let $\a\in W^{1/2,2}(\partial\Omega)$ have zero
total flux:
\begin{equation}
\label{fl0} \int\limits_{\partial\Omega}{\bf a}\cdot{\bf n}\,ds=0.
\end{equation}
Take an extension $\A$ satisfying
\begin{equation}
\label{AAA}
\begin{array}{r@{}l}
\A&{}\in W^{1,2}(\Omega),\\[2pt]
\div\A   &{}= 0\,\quad\hbox{\rm in } \Omega,\\[2pt]
\A &{}={\bf a} \quad\hbox{\rm on }
\partial\Omega,
\\[2pt]
\A(x) &{}={\bf 0} \quad\hbox{\rm if } x\in\R^2\setminus B_{R_0},
\end{array}
\end{equation}
where $B_{R_0}=B(0,R_0)$ is a disk of sufficiently large radius
such that
$$\frac12 B_{R_0}\supset\partial\Omega$$
(such extension exists because of condition~(\ref{fl0}), see,
e.g., \cite{LadSol}\,).

By a {\it weak solution} \,(\,$=D$-solution\,) of
problem~(\ref{SNS}) we mean a function ${\bf u}$ such that
$\u=\w+\A$, \ $\w\in H(\Omega)$, and the integral identity
\begin{equation}\label{weaksolution}
\nu\int\limits_\Omega\nabla{\bf u}\cdot\nabla\bfeta
\,dx+\int\limits_\Omega\big({\bf u}\cdot \nabla\big){\bf
u}\cdot\bfeta\,dx=0
\end{equation}
holds for any $\bfeta\in J_0^\infty(\Omega)$, where $
J_0^\infty(\Omega)$ is a set of all infinitely smooth solenoidal
vector-fields with compact support in $\Omega$. In particular, by
this definition we have
\begin{equation}\label{dd1}
\int\limits_{\Omega}|\nabla\u|^2dx<\infty.
\end{equation}  Moreover, by classical
regularity results for the Navier--Stokes system (see, e.g.,
\cite{OAL}, \cite{Galdibook}\,) every such solution is
$C^\infty$--regular inside the domain.

We look for a~solution to~(\ref{SNS}) as a limit of weak solutions
to the Navier--Stokes problem in a sequence of bounded domain
$\Omega_{bk}$ that in the limit exhaust  the unbounded domain
$\Omega$. The following result concerning the solvability of the
Navier-Stokes problem in bounded multi connected domains was
proved in~\cite{KPRAM}.

\begin{Theorem} \label{Th_Ex_b} {\sl
Let $\Omega'=\Omega_0\setminus
\bigl(\bigcup\limits_{j=1}^N\overline\Omega_j\bigr)$ be a~bounded
domain in $\R^2$ with multiply connected $C^2$-smooth boundary
$\partial\Omega'$ consisting of $N+1$ disjoint components
$\Gamma_j=\partial\Omega_j$, $j=0,\dots, N$. If~${\bf a}\in
W^{1/2,2}(\partial\Omega')$ satisfies
$$
\int\limits_{\partial\Omega'}{\bf a}\cdot{\bf n}\,ds=0,
$$
then $(\ref{SNS})$ with $\Omega=\Omega'$ admits at least one weak
solution ${\bf u}\in W^{1,2}(\Omega')$. }
\end{Theorem}

\begin{Remark}\label{ex-rem1}{\rm
Formally in the formulation of the existence theorem
in~\cite{KPRAM} we assumed that the boundary value~$\a$ satisfies
$\a\in W^{3/2,2}(\Omega)$ in order to have the regularity
condition~$(\u,p)\in W^{2,2}(\Omega)$. But really we used only
{\it local} variant of such regularity $(\u,p)\in
W^{2,2}_\loc(\Omega)$ (see~\cite[page~784, line 8 from
below]{KPRAM}\,). Now in our situation every $D$-solution has much
better $C^\infty$ regularity inside the domain~$\Omega$, so we
could assume less restrictive condition~$\a\in W^{1/2,2}(\Omega)$.}
\end{Remark}

Consider the sequence of boundary value problems
\begin{equation}
\label{NSESS} \left\{\begin{array}{r@{}l} - \nu \Delta\widehat{\bf
u}_k + (\widehat{\bf u}_k \cdot\nabla)\widehat{\bf u}_k +\nabla
\widehat p_k &{}= {\bf 0} \quad\hbox{\rm in }
\Omega_{bk},\\[2pt]
\div\widehat{\bf u}_k   &{}= 0\,\quad\hbox{\rm in } \Omega_{bk},\\[2pt]
\widehat{\bf u}_k &{}={\bf a} \quad\hbox{\rm on }
\partial\Omega,
\\[2pt]
\widehat{\bf u}_k &{}={\bf 0} \quad\hbox{\rm on }
\partial B_k=C_k,
\end{array}\right.
\end{equation}
where $\Omega_{bk}=B_k\cap\Omega$ for $k\ge k_0$, $B_k=\{x:\ \
|x|< k\}$, $\frac12B_{k_0}\supset
\bigcup\limits_{i=1}^N\overline\Omega_i$. By
Theorem~\ref{Th_Ex_b}, each problem~(\ref{NSESS}) has a~solution
$\widehat{\bf u}_k\in W^{1,2}(\Omega_{bk})$ satisfying
$\div\widehat{\bf u}_k=0$ and the corresponding integral
identities (of \eqref{weaksolution} type).

Assume that there is a positive constant $c$ independent of $k$
such that
\begin{equation}
\label{kkkaf} \int\limits_\Omega|\nabla\widehat{\bf u}_k|^2dx\le c
\end{equation}
(possibly along a subsequence of $\{\widehat{\bf
u}_k\}_{k\in{\mathbb N}}$). This estimate
 implies the existence of a
solution to problem~(\ref{SNS}). Indeed, from~(\ref{kkkaf}) and
from the boundary conditions~(\ref{NSESS}${}_3$) it follows that
the~sequence $\widehat{\bf u}_k$ is bounded in
$W^{1,2}_\loc(\overline{\Omega})$. Hence, $\widehat{\bf u}_k$
converges weakly (modulo a subsequence) in
$W^{1,2}_\loc(\overline{\Omega})$ and strongly in $L^q_{\rm
loc}(\overline\Omega)$ $(1\le q<\infty)$ to a function
$\widehat{\bf u}\in D^{1,2}_\sigma(\Omega)$.  It is easy to check
that this limiting function~$\widehat\u$ is a~$D$-solution to the
 Navier--Stokes problem~(\ref{SNS})
in the exterior domain $\Omega$.

 Thus, to prove the
assertion of Theorem~\ref{T2}, it is sufficient to establish the
uniform estimate~(\ref{kkkaf}). We shall prove (\ref{kkkaf})
following a classical {\it reductio ad absurdum\/} argument of
J.~Leray \cite{Ler} and O.A.~Ladyzhenskaia \cite{OAL}. If
(\ref{kkkaf}) is not true, then there exists a sequence
$\{\widehat\u_k\}_{k\in{\mathbb N}}$ such that
$$
\lim_{k\to+\infty}J_k^2=+\infty,\quad
J_k^2=\int\limits_\Omega|\nabla\widehat{\bf u}_k|^2dx.
$$
The sequence ${\bf u}_k=\widehat{\bf u}_k/J_k$ is   bounded in
$D^{1,2}_\sigma(\Omega)\cap L^q_\loc(\overline\Omega)$ and it
holds
\begin{equation}
\label{NSESS'}
 \dfrac{\nu}{J_k}\int\limits_{\Omega}\nabla{\bf u}_k\cdot\nabla\bfeta\,dx= -\int\limits_\Omega ({\bf u}_k\cdot\nabla){\bf u}_k\cdot\bfeta\,dx
 \end{equation}
for all $\bfeta\in H(\Omega_{bk})$. Extracting  a subsequence (if
necessary) we can assume that ${\bf u}_k$ converges weakly in
$D^{1,2}_\sigma(\Omega)$ and strongly in $L^q_{\rm loc}
(\overline\Omega)$ $(1\le q<\infty)$  to a vector field ${\bf
v}\in H(\Omega)$ with
\begin{equation}\label{norm1}
\int\limits_\Omega|\nabla {\bf v}|^2dx\le1.
\end{equation} Fixing in \eqref{NSESS'} a
solenoidal smooth  $\bfeta$ with compact support and letting
$k\to+\infty$ we get
\begin{equation}\label{Eweak}
\int\limits_\Omega({\bf v}\cdot\nabla){\bf
v}\cdot\bfeta\,dx=0\quad\forall \,\bfeta\in J_0^\infty(\Omega),
\end{equation}
Hence, ${\bf v}\in H(\Omega)$ is a weak solution to the Euler
equations, and for some $ p\in W^{1,q}_\loc(\overline\Omega)$, \
($1<q<\infty$), \ the pair $({\bf v},  p)$ satisfies the Euler
equations almost everywhere:
\begin{equation}
\label{2.1}\left\{\begin{array}{rcl} \big({\bf
v}\cdot\nabla\big){\bf v}+\nabla p & = & 0 \qquad \ \ \
\hbox{\rm in }\;\;\Omega,\\[4pt]
\div{\bf v} & = & 0\qquad \ \ \ \hbox{\rm in }\;\;\Omega,
\\[4pt]
{\bf v} &  = & 0\ \
 \qquad\  \hbox{\rm on }\;\;\partial\Omega.
\end{array}\right.
\end{equation}

Put $\nu_k=(J_k)^{-1}\nu$. Then the system~(\ref{NSESS})  could be
rewritten in the following form

\begin{equation}
\label{NSk} \left\{\begin{array}{r@{}l} - \nu_k \Delta{\bf u}_k +
({\bf u}_k \cdot\nabla){\bf u}_k +\nabla  p_k &{}= {\bf 0}
\quad\hbox{\rm in }
\Omega_{bk},\\[2pt]
\div{\bf u}_k   &{}= 0\,\quad\hbox{\rm in } \Omega_{bk},\\[2pt]
{\bf u}_k &{}=\frac{\nu_k}{\nu}{\bf a} \quad\hbox{\rm on }
\partial\Omega,
\\[2pt]
{\bf u}_k &{}={\bf 0} \quad\hbox{\rm on }
\partial B_k=C_k.
\end{array}\right.
\end{equation}
where $\ue_k,\,p_k\in C^\infty_{\loc}(\Omega_{bk})$. In conclusion,  we
come to the following assertion.

\begin{Lemma}
\label{lem_Leray_symm} {\sl Assume that  $\Omega\subset\R^2$ is an
exterior   domain of type \eqref{Omega} with $C^2$-smooth boundary
$\partial\Omega$,  and ${\bf a}\in W^{1/2,2}(\partial\Omega)$
satisfies  zero total flux condition~(\ref{dom5}). If the
assertion of Theorem~\ref{T2} is false, then there exist $ \ve, p$
with the following properties.

\medskip

(E) \ \,The functions $\ve\in H(\Omega)$, \ $p\in
W^{1,q}_\loc(\overline\Omega)$, \  ($1<q<\infty$) satisfy the
Euler system \eqref{2.1}.

\medskip

(E-NS) \ \,Condition (E) is fulfilled and there exist a sequences
of functions $\ue_k \in W^{1,2}(\Omega_{bk})$, $p_k\in
{W^{1,q}(\Omega_{bk})}$,  $\Omega_{bk}=\Omega\cap B_{R_k}$,
$R_k\to\infty$ as $k\to\infty$, and numbers $\nu_k\to0+$, such
that the pair $(\ue_k,p_k)$ satisfies~(\ref{NSk}), and
\begin{equation}
\label{E-NS-ax} \|\nabla\ue_k\|_{L^2(\Omega_{bk})}\equiv1,\quad
\ue_k\rightharpoonup \ve\mbox{ \ in \
}W_\loc^{1,2}(\overline\Omega),\quad p_k\rightharpoonup p\mbox{ \
in \ }W_\loc^{1,q}(\overline\Omega).
\end{equation}
\begin{equation}\label{cont_e}
\begin{array}{l}
\displaystyle \nu= \int\limits_\Omega ({\bf v}\cdot\nabla){\bf
v}\cdot\A\,dx
\end{array}\end{equation}
Moreover, $\ue_k,\,p_k\in C^\infty(\Omega_{bk})$ {\rm(}this notation
means $C^\infty$-regularity inside the domain~$\Omega_{bk}$\,{\rm)}.}
\end{Lemma}

\proof

We need to prove only the identity~(\ref{cont_e}), all other
properties are already established above. By construction
$\u_k=\w_k+\frac1{J_k}\A$, where $\w_k\in H(\Omega_{bk})$, in
particular, $\w_k\equiv0$ on $\partial\Omega_{bk}$.
 Choosing $\bfeta={\bf
w}_k$ in (\ref{NSESS'}) and integration by parts yields
 \begin{equation}
\label{NaaS} \displaystyle \nu=\int\limits_\Omega ({\bf
w}_k\cdot\nabla){\bf w}_k\cdot{\bf A} \,dx+ {1\over
J_k}\int\limits_\Omega{\bf A}\cdot\nabla{\bf w}_k\cdot{\bf A} \,dx
+{\nu\over J_k}\int\limits_\Omega\nabla{\bf A} \cdot\nabla{\bf u}_k\,dx.
\end{equation}
Since $\A\in W^{1,2}(\Omega)$  has a compact support, it is easy
to check that we can pass to the limit in~(\ref{NaaS}) and receive
the~required assertion~(\ref{cont_e}). $\qed$

\medskip

Notice that because of~(\ref{cont_e}) the limiting solution~$\ve$
of the Euler system \eqref{2.1} is nontrivial.

Now, to finish the proof of Theorem~\ref{T2}, we need to show that
conditions~(E-NS) lead to a contradiction. The next two
subsections are devoted to this purpose.

\subsection{Some properties of solutions to Euler system}
\label{poet2}

In this section we assume that the assumptions~(E) of
Lemma~\ref{lem_Leray_symm} are satisfied.  In particular,

\begin{equation}
\label{ax'3}\int\limits_{{\Omega}}|\nabla{\bf v}(x)|^2\,dx<\infty.
\end{equation}

The next statement  was proved in \cite[Lemma 4]{KaPi1} and in
\cite[Theorem 2.2]{Amick}.

\begin{Theorem}
\label{kmpTh2.3'} {\sl Let the conditions {\rm (E)} be fulfilled.
Then
\begin{equation} \label{bp2} \forall j\in\{1,\dots,N\} \ \exists\, \widehat p_j\in\R:\quad
p(x)\equiv \widehat p_j\quad\mbox{for }\Ha^1-\mbox{almost all }
x\in\Gamma_j.\end{equation} }
\end{Theorem}

Using the last fact, below we assume without loss of generality
that the functions $\ve,p$ are extended to the whole plane $\R^2$
as follows:
\begin{equation}
\label{axc10.10} \ve(x):=0,\quad x\in \R^2\setminus\Omega,
\end{equation}
\begin{equation}
\label{axc110} p(x):=\widehat p_j,\ \, x\in \R^2\cap\bar\Omega_j,\
j=1,\dots,N.
\end{equation}
Obviously,  the extended functions inherit the properties of the
previous ones. Namely, $\ve\in H(\R^2)$, $p\in
W^{1,q}_\loc(\R^2)$, and the Euler equations~(\ref{2.1}) are
fulfilled almost everywhere in~$\R^2$. That means, the pair
$(\ve,p)$ is a weak (=Sobolev) solution to Euler
system~(\ref{2.1}) {\it in the whole plane}.

First of all all, we prove the uniform boundedness and continuity
of  the~pressure.

 \begin{Theorem}
\label{pp} {\sl Let the conditions {\rm (E)} be fulfilled. Then
\begin{equation} \label{co-1}p\in
D^{2,1}(\R^2)\cap D^{1,2}(\R^2).
\end{equation}
In particular, the function~$p$ is continuous and convergent at
infinity, i.e.,
\begin{equation} \label{co-2}
\exists \lim\limits_{x\to\infty}p(x)\in \R.
\end{equation}}
\end{Theorem}

\proof   By well-known fact concerning $D$-solutions to Euler and
Navier--Stokes system (see, e.g.,\cite[Lemma~4.1]{GW2}), the
averages of the pressure are uniformly bounded:
\begin{equation} \label{co-3}
\sup\limits_{r>0} \biggl|\frac1r\int\limits_{C_r}p\,ds\biggr|<\infty,
\end{equation}
where, recall, $C_r=\{x\in\R^2:|x|=r\}$. Moreover, since
$\int\limits_{\R^2}|\nabla\ve|^2dx<\infty$, there exists an
increasing sequence $r_i\to +\infty$ such that
\begin{equation} \label{co-4}
\int\limits_{C_{r_i}}|\nabla\ve|\,ds\le \e_i\to0\quad\mbox{ as
}i\to\infty
\end{equation}
and
\begin{equation} \label{co-4'}
\sup\limits_{x\in C_{r_i}} \bigl|\ve(x)\bigr|\le \e_i\sqrt{\ln
r_i}
\end{equation}
(see \cite[Lemmas~2.1--2.2]{GW2})). From
(\ref{co-3})--(\ref{co-4}) and from the equation~(\ref{2.1}$_1$)
it follows that
\begin{equation} \label{co-5}
\sup\limits_{x\in C_{r_i}} \bigl|p(x)\bigr|\le C\,\sqrt{\ln r_i}.
\end{equation}
Indeed,
$$
|p(r_i,\theta)-\bar p(r_i)|\le \int\limits_{C_{r_i}}|\nabla p|\,ds
\leq \int\limits_{C_{r_i}}|\ve|\cdot|\nabla \ve|\,ds\leq \e_i\sqrt{\ln
r_i} \int\limits_{C_{r_i}}|\nabla \ve|\,ds\le \e^2_i\sqrt{\ln r_i},
$$
here $\bar p(r_i)=\frac{1}{2\pi r_i}\int\limits_{C_{r_i}}p\,ds$.  The last
inequality and the uniform boundedness of $\bar p(r_i)$
(see~(\ref{co-3})\,) implies~\eqref{co-5}.

Clearly, $p\in W^{1,q}_\loc(\R^2)$ is the weak solution to the
Poisson equation
\begin{equation} \label{r2.2}
\Delta p =-\nabla{\mathbf v}\cdot\nabla{\mathbf
v}^\top\quad\mbox{in }\R^2
\end{equation}
(recall that after our agreement about extension of ${\mathbf v}$
and $p$, see (\ref{axc10.10})--(\ref{axc110}), the Euler equations
(\ref{2.1}) are fulfilled in the whole $\R^2$).

Put
$$
G(x)=-{1\over 2\pi}\int\limits_{\Omega} \log|x-y|(\nabla{\mathbf
v}\cdot\nabla{\mathbf v}^\top)(y) dy.
$$

 By the results of \cite{CLMS}, $\nabla{\mathbf
v}\cdot\nabla{\mathbf v}^\top$  belongs to the Hardy space
$\H^1(\R^2)$. Hence by Calder\'on--Zygmund theorem for Hardy's
spaces \cite{Stein} $G\in D^{2,1}(\R^2)\cap D^{1,2}(\R^2)$. By
classical facts from the theory of Sobolev Spaces (see, e.g.,
\cite{maz'ya}\,), the last inclusion implies that~$G$ is
continuous and convergent at infinity, in particular,
\begin{equation} \label{co-6}
\sup\limits_{x\in\R^2}|G(x)|<\infty.
\end{equation}
Consider the function $p_*=p-G$. By construction, $\Delta p_*=0$
in $\R^2$, i.e., $p_*$ is a harmonic function, and
by~(\ref{co-5}), (\ref{co-6}) we have
\begin{equation} \label{co-7}
\sup\limits_{x\in C_{r_i}}|p_*(x)|\le C\sqrt{\ln r_i}.
\end{equation}
 From the Liouville type theorems for harmonic functions  (see, i.e.,
\cite{Axler}\,) it follows that $p_*\equiv\const$. Consequently,
$p\equiv G+\const$, that implies the assertions of the Theorem.
$\qed$

\bigskip

We say that the function $f\in W^{1,s}_\loc(\R^2)$ satisfies a
{\it weak one-side maximum principle}, if
\begin{equation}
\label{2.13} \mathop{\hbox{\rm ess}\,\hbox{\rm
sup}}_{x\in\Omega^\prime}\,f(x)\leq \mathop{\hbox{\rm
ess}\,\hbox{\rm sup}}_{x\in\partial\Omega^\prime}\,f(x)
\end{equation}
holds for any bounded subdomain $\Omega^\prime$ with the boundary
$\partial\Omega^\prime$ not containing singleton connected
components.  (In (\ref{2.13}) negligible sets are the sets of
2--dimensional Lebesgue measure zero in the left \emph{esssup},
and the sets of 1--dimensional Hausdorff measure zero in the right
\emph{esssup}.)

The total head pressure for the Euler system
$$
\Phi:=p+\frac12|\ve|^2.
$$
plays an important role in the forthcoming considerations. The
following two results were proved in~\cite{kpr}.

\begin{Theorem}
\label{mp1}{\sl Suppose that the assumptions~$(E-NS)$ from the
previous subsection are satisfied. Than the total head pressure
$\Phi$ satisfies the weak maximum principle in $\R^2$.}
\end{Theorem}

The second equality in (\ref{2.1}) (which is fulfilled, after the
above extension agreement, see~(\ref{axc10.10})--(\ref{axc110}),
in the whole plane~$\R^2$) implies the existence of a stream
function $\psi\in W^{2,2}_{\loc}(\R^2)$ such that
\begin{equation}
\label{ax7}\nabla\psi=\ve^\bot,
\end{equation}
i.e.,
\begin{equation}
\label{ax7'}\frac{\partial\psi}{\partial x_1}=v_2,\qquad
\frac{\partial\psi}{\partial x_2}=-v_1.
\end{equation}

Let us formulate regularity results concerning the considered
functions.

\begin{Lemma}[see, e.g., Theorem 3.1 in \cite{kpr}]\label{kmpTh2.1}
{\sl If conditions {\rm (E)} are satisfied, then $\psi\in C(\R^2)$
and there exists a set $A_{\ve}\subset \R^2$ such that

 {\rm (i)}\quad $ \mathfrak{H}^1(A_{\ve})=0$;

{\rm (ii)} for  all  $x\in\Omega\setminus A_{\ve}$
\begin{displaymath}
\lim\limits_{r\to
0}\dashint\nolimits_{B_r(x)}|\ve(z)-\ve(x)|^2dz=\lim\limits_{r\to
0}\dashint\nolimits_{B_r(x)}|{\Phi}(z)-{\Phi}(x)|^2dz=0;
\end{displaymath}
moreover, the function $\psi$ is differentiable at $x$ and
$\nabla\psi(x)=(v_2(x),- v_1(x))$;

{\rm  (iii) } for every  $\varepsilon >0$ there exists a set
$U\subset \mathbb{R}^2$ with
$\mathfrak{H}^1_\infty(U)<\varepsilon$ such that $A_{\ve}\subset
U$ and the functions $\ve, \Phi$ are continuous in $\R^2\setminus
U$.}
\end{Lemma}

By virtue of~(\ref{axc10.10}), we have $\nabla \psi(x)=0$ for
almost all $x\in{\Omega}_j$. Then
\begin{equation}
\label{axc10} \forall j\in\{1,\dots,N\}\ \ \exists\, \xi_j\in\R:
 \quad \psi(x)\equiv\xi_j\qquad \forall x\in \overline{\Omega_j}\cap\R^2.
\end{equation}

By direct calculations one easily gets the identity
\begin{equation}
\label{2.2} \nabla\Phi=\omega\nabla\psi,
\end{equation}
here $\omega=\Delta\psi=\partial_1 v_2-\partial_2 v_1$ means the
corresponding vorticity.

The next assertion, obtained in the paper~\cite{kpr}, is the
another important tool  for the proof of Theorem 2.

\medskip

\begin{Theorem}[{\bf Bernoulli Law for Sobolev solutions}]
\label{kmpTh2.2} {\sl Let the conditions~{\rm (E)} be valid. Then
there exists a set $A_\ve\subset \R^2$ with $\Ha^1(A_\ve)=0$, such
that for any compact connected\footnote{We understand the
connectedness in the sense of general topology.} set
$K\subset\R^2$ the following property holds : if
\begin{equation}
\label{2.4} \psi\big|_{K}=\const,
\end{equation}
then
\begin{equation}
\label{2.5'}  \Phi(x_1)=\Phi(x_2) \quad\mbox{for
 all \,}x_1,x_2\in K\setminus A_{\bf v}.
\end{equation}
}
\end{Theorem}

Of course, we could assume without loss of generality that the
sets $A_\ve$ from Lemma~\ref{kmpTh2.1} and Theorem~\ref{kmpTh2.2}
are the same.

Identities (\ref{axc10.10})--(\ref{axc110}) mean that
\begin{equation}
\label{axc11} \Phi(x)\equiv\widehat p_j\qquad \forall x\in
\R^2\cap\overline{\Omega_j},\ \,j=1,\dots,N.
\end{equation}

By Theorem~\ref{pp} (see also its proof) there exists an
increasing  sequence of numbers $r_i\to+\infty$ such that
\begin{equation}
\label{axc17} \sup\limits_{x\in
C_{r_i}}|\Phi(x)-\overline{\Phi_i}| \to 0,
\end{equation}
where $\overline{\Phi_i}=\dashint_{C_{r_i}}\Phi(x)\,ds$ is the mean
value of $\Phi$ over the circle~$C_{r_i}$. Indeed, by definition,
$|\nabla\Phi|\le |\ve|\cdot|\nabla\ve|$. By standard estimates
(e.g., \cite[Lemma~2.1]{Galdibook}\,)
\begin{equation}
\label{axc17-check1}\int\limits_{C_r}|\ve|^2ds\le Cr\ln{r}.
\end{equation}
Further, since
$\int\limits_{\R^2}|\nabla\ve|^2dx=\int\limits_0^\infty
dr\,\int\limits_{C_r}|\nabla\ve|^2ds<\infty$, there exists an
increasing sequence $r_i\to +\infty$ such that
\begin{equation}
\label{axc17-check2}\int\limits_{C_{r_i}}|\nabla\ve|^2ds\le
\frac{\e_i}{r_i\ln{r_i}},
\end{equation}
where $\e_i\to0$. Formulas
(\ref{axc17-check1})--(\ref{axc17-check2}) and the H\"older
inequality imply
\begin{equation}
\label{axc17-check3}\int\limits_{C_{r_i}}|\nabla\Phi|\,ds\le
\int\limits_{C_{r_i}}|\ve|\cdot|\nabla\ve|\,ds\le\sqrt{C\e_i}\to 0,
\end{equation}
thus we obtain~(\ref{axc17}).

From the weak maximum principle (see Theorem~\ref{mp1}\,) it
follows that there exists a limit
$\Phi_\infty=\lim\limits_{i\to\infty}\overline{\Phi_i}$, which
does not depend on the choice of circles~$C_{r_i}$ (it can
be~$\Phi_\infty=\infty$). Again the same maximum principle
implies that
\begin{equation}
\label{axc18} \esssup\limits_{x\in
\R^2}\Phi(x)=\max\{\Phi_\infty,\widehat p_1,\dots,\widehat
p_N\},\end{equation} where $\widehat p_j$ are the constants form
Theorem~\ref{kmpTh2.3'}. Further we consider separately three
possible cases.

(a) The maximum of $\Phi$ is attained strictly at
infinity\footnote{The case
$\esssup\limits_{x\in\Omega}\Phi(x)=+\infty$ is not excluded. },
i.e.,
\begin{equation}\label{as-prev1}
\Phi_\infty=\esssup\limits_{x\in\Omega}\Phi(x)>\max\{\widehat
p_1,\dots,\widehat p_N\}.
\end{equation}

(b) The maximum of $\Phi$ is attained on some boundary
component~--- not at infinity:
\begin{equation}\label{as1-axxx} \max\{\widehat p_1,\dots,\widehat p_N\}=\esssup\limits_{x\in\Omega}\Phi(x)>\Phi_\infty.
\end{equation}

(c) The maximum of $\Phi$ is attained both at infinity and on some
boundary component:
\begin{equation}\label{as-prev-id-ax}
\Phi_\infty=\esssup\limits_{x\in\Omega}\Phi(x)=\max\{\widehat
p_1,\dots,\widehat p_N\}.
\end{equation}

\subsection{The case $\esssup\limits_{x\in\Omega}\Phi(x)=\Phi_\infty>\max\{\widehat p_1,\dots,\widehat p_N\}$.}
\label{EPcontr-axx}

Let us consider  the first case~(\ref{as-prev1}).

We will adopt the arguments of~\cite[subsection 2.4.1]{KPRAM}.
Note that the calculation in the present situations are much
easier, since the set where $\Phi$ close to the maximum is
separated from the boundary components. For the reader convenience, in this subsection we reproduce these arguments
in details.

Without loss of generality we could assume that
\begin{equation}\label{cas1}
\Phi_\infty>\delta>0>-\delta>\max\{\widehat p_1,\dots,\widehat
p_N\},
\end{equation}
where $\delta$ is sufficiently small positive number.

By definition of~$\Phi_\infty$ (see, e.g., (\ref{axc17})\,), there
exists a radius $r_0>0$ such that~$B_{\frac12
r_0}\supset\partial\Omega$ and

\begin{equation}\label{cas2}
C_{r_0}\cap A_\ve=\emptyset;
\end{equation}

\begin{equation}\label{cas3}
\inf\limits_{x\in C_{r_0}}\Phi(x)\ge\delta.
\end{equation}

Our first goal is to separate the boundary components $\Gamma_j$
where $\Phi<0$ from $C_{r_0}$ by level sets of~$\Phi$ compactly
supported in $\Omega$. More precisely, for any $t\in(0,\delta)$
and $j=1,\dots,N$ we construct a~continuum $A_j(t)\Subset \Omega$
with the following properties:

(i) The set $\Gamma_j=\partial\Omega_j$ lies in a~bounded
connected component of the open set~$\R^2\setminus A_j(t)$;

(ii) $\psi|_{A_j(t)}\equiv\const$, \ $\Phi(A_j(t))=-t$;

(iii) (monotonicity) If $0<t_1<t_2<\delta_{p}$, then $A_j(t_1)$
lies in the unbounded connected component of the
set~$\R^2\setminus A_j(t_2)$ (in other words, the set
$A_j(t_2)\cup\Gamma_j$ lies in the bounded connected component of
the set~$\R^2\setminus A_j(t_1)$, see Fig.1).
\begin{center}
\includegraphics[scale=0.4]{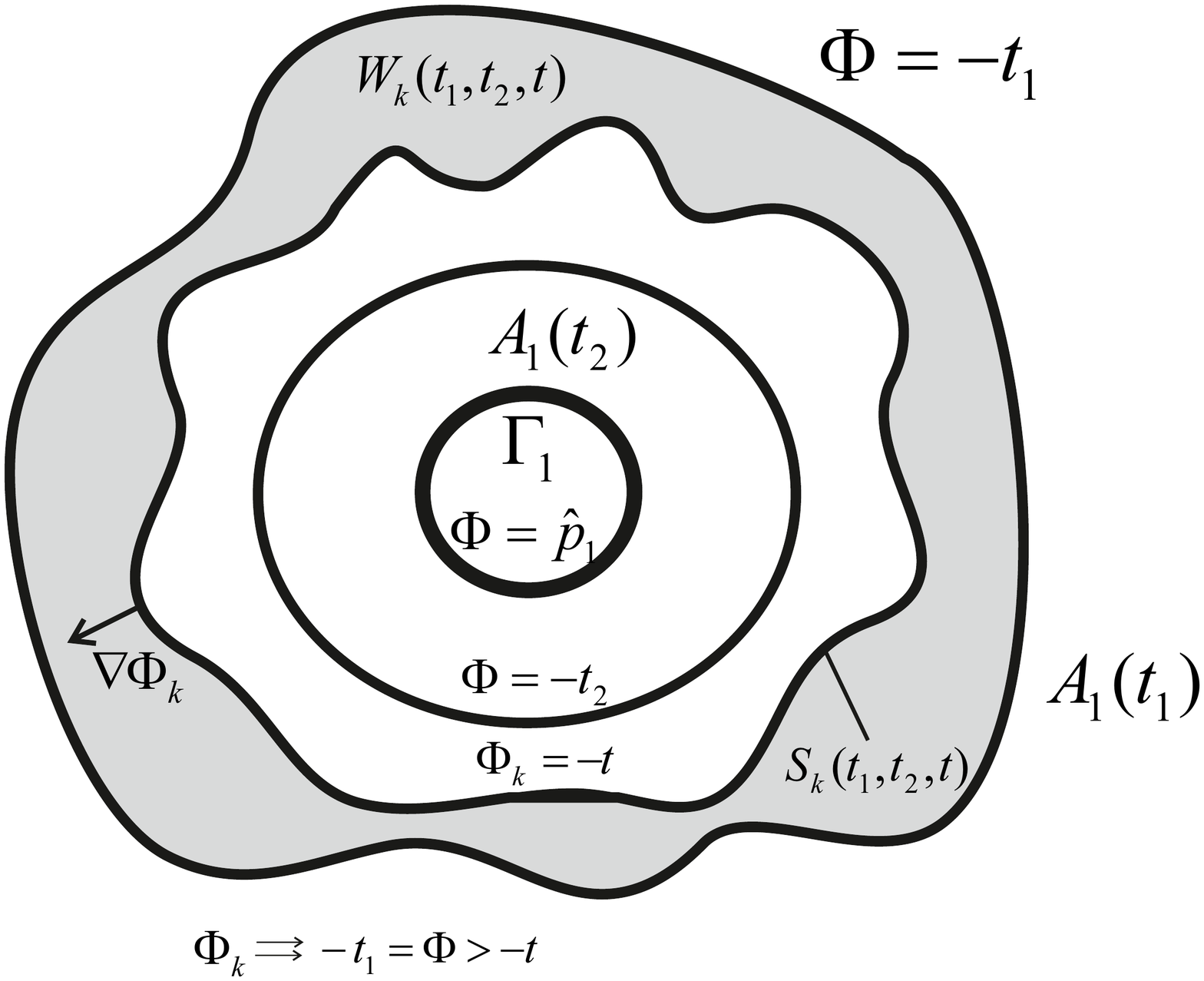}
\end{center}
\begin{center}
Fig. 1. {\sl The surface $S_k(t_1,t_2,t)$ for the case of $N=1$.}
\end{center}

For this construction, we shall use the results of
Subsection~\ref{Kronrod-s}. More precisely, we apply Kronrod's
results to the stream function~$\psi|_{\bar B_{r_0}}$.
Accordingly, $T^0_{\psi}$ means the corresponding Kronrod tree for
the restriction~$\psi|_{\bar B_{r_0}}$.

For any element~$C\in T^0_\psi$ with $C\setminus
A_\ve\ne\emptyset$  we can define the value $\Phi(C)$ as
$\Phi(C)=\Phi(x)$, where $x\in C\setminus A_\ve$. This definition
is correct because of the Bernoulli Law. (In particular, $\Phi(C)$
is well defined if $\diam C>0$.)

Take points $x_0\in C_{r_0}$ and $x_j\in \Omega_j$, $j=1,\dots,N$,
such that the straight segment $L_j$ with endpoints $x_0$ and
$x_j$ satisfies
\begin{equation}\label{f-m2}
L_j\cap A_\ve=\emptyset;\end{equation}
\begin{equation}\label{cas5}
\mbox{the restriction $\Phi|_{L_j}$ is a continuous
function}\end{equation} (the existence of such points and segments
follows from Lemma~\ref{kmpTh2.1}~(iii)\,).

Denote by $E_{0}$ and $E_j$ the elements of $T^0_\psi$ with
$x_0\in E_0$ and $x_j\in E_j$.  Note that from
$\psi|_{\Omega_j}\equiv\const$ it follows that $\overline\Omega_j\subset
E_j$.  Consider the arc~$[E_j,E_{0}]\subset T^0_\psi$. Recall
that, by definition, a connected component $C$ of a level set
of~$\psi|_{\bar B_{r_0}}$ belongs to the arc~$[E_j,E_{0}]$ iff
$C=E_{0}$, or $C=E_j$, or $C$ separates $E_{0}$ from $E_j$ in
$\bar B_{r_0}$, i.e., if $E_0$ and $E_j$ lie in different
connected components of~$\bar B_{r_0}\setminus C$. In particular,
since $E_0\cap L_{j}\ne\emptyset\ne E_j\cap L_{j}$, we have
\begin{equation}\label{f-m1}
C\cap L_{j}\ne\emptyset\quad\forall C\in [E_j,E_0].\end{equation}
Therefore, in view of equality~(\ref{f-m2}) the value $\Phi(C)$ is
well defined for all~$C\in [E_j,E_0]$. Moreover, we have
\begin{Lemma}
\label{l-m1} {\sl The restriction $\Phi|_{[E_j,E_0]}$ is a
continuous function. }
\end{Lemma}

\proof The assertion follows immediately\footnote{See also the
proof of Lemma 3.5 in \cite{KPRAM}.} from the assumptions
(\ref{f-m2})--(\ref{f-m1}), from the continuity of $\Phi|_{L_j}$,
and from the definition of convergence in~$T^0_\psi$ (see
Subsection~\ref{Kronrod-s}\,). \qed

Define the natural order\footnote{Recall, that by
Lemma~\ref{kmpLem6}, the set $[E_j,E_0]$  is homeomorphic to the
segment of a real line, i.e. it is an arc. So we could define
a~natural order on this arc and take maxima, minima etc.~--- as
for usual segment. There are two symmetric possibilities to define
a usual linear order on the arc; here by our choice $E_j<E_0$.} on
the arc~$[E_j,E_0]$. Namely, we say, that $A< C$ for some
different elements~$A,C\in[E_j,E_0]$ iff $C$ closer to~$E_0$ than
$A$, i.e., if the sets $E_0$ and $C$ lie in the same connected
component of the set~$\bar B_{r_0}\setminus A$.

Put $$K_{j}=\min\{K\in [E_j,E_0]: K \cap
 C_{r_0}\ne\emptyset\}$$
 (this minimum exists since $E_0\cap C_{r_0}\ne\emptyset$\,).
By elementary and obvious topological arguments we have
\begin{equation}\label{casf1}
\forall K\in [E_j,E_0]\quad\bigl(\,K \cap
 C_{r_0}\ne\emptyset\Leftrightarrow K\ge K_j\bigr).\end{equation}
From (\ref{cas2})--(\ref{cas3}) and from the Bernoulli Law it follows that
\begin{equation}\label{casf2}
\Phi(K)\ge\delta\qquad\forall K\in [K_j,E_0].\end{equation} In
particular, since $\Phi(E_j)<-\delta$, we have
\begin{equation}\label{casf3}
E_j<K_j\le E_0.\end{equation} By construction,
\begin{equation}\label{casf4}
K\cap C_{r_0}=\emptyset\qquad\forall K\in [E_j,K_j),\end{equation}
where, as usual, $[E_j,K_j)=[E_j,K_j]\setminus\{K_{j}\}$.

We say that a set $\mathcal Z\subset [E_j,E_0]$ has $T$-measure
zero if $\Ha^1(\{\psi(K):K\in \mathcal Z\})=0$.

\begin{Lemma}
\label{regc-ax} {\sl For every $j=1,\dots,N$,  $T$-almost all
$K\in[E_j,K_j]$ are $C^1$-curves homeomorphic to the circle and
$K\cap A_\ve=\emptyset$. Moreover, there exists
a~subsequence~$\Phi_{k_l}$ such that $\Phi_{k_l}|_K$ converges to
$\Phi|_K$ uniformly $\Phi_{k_l}|_K\rightrightarrows\Phi|_K$ on
$T$-almost all $K\in[E_j,E_0]$.}
\end{Lemma}

\proof The first assertion of the lemma follows from
Theorem~\ref{kmpTh1.1}~(iii) and ~(\ref{casf4}). The validity of
the second one for $T$-almost all
$K\in[E_j,K_j]$ was proved in~\cite[Lemma 3.3]{kpr}. \qed \\

Below we assume (without loss of generality) that the subsequence
$\Phi_{k_l}$ coincides with the whole sequence $\Phi_{k}$.
Furthermore, we will call {\it regular } the cycles $K$ which
satisfy the~assertion of Lemma~\ref{regc-ax}.

Since $\diam C>0$ for every $C\in [E_j,E_0]$, we obtain, by
\cite[Lemma~3.6]{KPRAM},  that the function
$\Phi|_{[E_j,E_0]}$ has the following analog of Luzin's
$N$-property.

\begin{Lemma}
\label{lkr7} {\sl For every $j=1,\dots,N$, if $\mathcal Z\subset
[E_j,E_0]$ has $T$-measure zero, then $\Ha^1(\{\Phi(K):K\in
\mathcal Z\})=0$.}
\end{Lemma}

Note that Lemma~\ref{lkr7} is not tautological: in the definition of
$T$-zero measure we have stream function~$\psi$, but
Lemma~\ref{lkr7} deals about another function, total head
pressure~$\Phi$. It looks like Luzin $N$-property: $\psi(E)$ has
zero measure implies $\Phi(E)$ has zero measure.

From Lemmas~\ref{regc-ax}--\ref{lkr7}  and from~(\ref{casf2}) we
conclude

\begin{Corollary}
\label{regPhi-ax} {\sl For every $j=1,\dots,N$ and for almost all
$t\in(0,\delta)$ we have
$$\bigl(\,K\in[E_j,E_0]\mbox{\rm\ and }\Phi(K)=-t\,\bigr)\Rightarrow
K\mbox{\rm\ is a regular cycle}.$$}
\end{Corollary}

Below we will say that a value~$t\in(0,\delta)$ is {\it regular}
if it satisfies the assertion of Corollary~\ref{regPhi-ax}. Denote
by~$\Ti$  the set of all regular values. Then the set
$(0,\delta)\setminus \Ti$ has zero measure.

For $t\in (0,\delta)$ and $j\in\{1,\dots,N\}$ denote
$$A_j(t)=\max\{K\in[E_j,E_0]:\Phi(K)=-t\}.$$
By construction, the function $A_j(t)$ is nonincreasing and
satisfies the properties (i)--(iii) from the beginning of this
subsection. Moreover, by definition of regular values we have the
following additional property:

(iv) If $t\in \Ti$, then $A_j(t)$ is a regular cycle\footnote{Some
of these cycles $A_j(t)$ could coincide, i.e., equalities of type
$A_{j_1}(t)=A_{j_2}(t)$ are possible (if Kronrod arcs
$[E_{j_1},E_0)$ and $[E_{j_2},E_0)$ have nontrivial intersection),
but this  a~priori possibility has no influence on our
arguments.}.

For $t\in\Ti$ denote by ${V}(t)$ the unbounded connected component
of the open set $\R^2\setminus\bigl(\cup_{j=1}^N A_j(t)\bigr)$.
Since $A_{j_1}(t)$ can not separate $A_{j_2}(t)$ from
infinity\footnote{Indeed, if $A_{j_2}(t)$ lies in a~bounded
component of $\R^2\setminus A_{j_1}(t)$, then by construction
$A_{j_1}(t)\in [E_{j_2},E_0]$ and $A_{j_1}(t)>A_{j_2}(t)$ with
respect to the~above defined order on~$[E_{j_2},E_0]$. However,  it
contradicts the definition of
$A_{j_2}(t)=\max\{K\in[E_{j_2},E_0]:\Phi(K)=-t\}$.} for
$A_{j_1}(t)\ne A_{j_2}(t)$, we have
\begin{equation}\label{boundary1}\partial{V}(t)=
A_{1}(t)\cup\dots\cup A_N(t),\quad\ t\in\Ti.\end{equation}  By
construction, the sequence of domains ${V}(t)$ is increasing,
i.e., ${V}(t_1)\subset{V}(t_2)$ for $t_1<t_2$.

Let $t_1,t_2\in\Ti$ and $t_1<t_2$. The next geometrical objects
plays an~important role in the estimates below: for
$t\in(t_1,t_2)$ we define the level set $S_k(t,t_1,t_2)\subset
\{x\in\Omega_{bk}:\Phi_k(x)=-t\}$ separating cycles $\cup_{j=1}^N
A_j(t_1)$ from $\cup_{j=1}^NA_j(t_2)$ as follows. Namely, take
arbitrary $t',t''\in\Ti$ such that $t_1<t'<t''<t_2$. From
Properties~(ii),(iv) we have the uniform convergence
$\Phi_k|_{A_j(t_1)}\rightrightarrows -t_1$,
$\Phi_k|_{A_j(t_2)}\rightrightarrows -t_2$ as $k\to\infty$ for
every~$j=1,\dots,N$. Thus there exists
$k_\circ=k_\circ(t_1,t_2,t',t'')\in\N$ such that for all $k\ge
k_\circ$
\begin{equation}\label{boundary0}\Phi_k|_{A_j(t_1)}>-t',\quad
\Phi_k|_{A_j(t_2)}< -t''\quad\forall j=1,\dots,N.\end{equation} In
particular,
\begin{equation}\label{boundary2}
\begin{array}{lcr}
\Phi_k|_{A_j(t_1)}> -t,\quad\Phi_k|_{A_j(t_2)}< -t,\quad \forall t\in[t',t''], \ \forall k\ge k_\circ,
\\
\\
\forall j=1,\dots,N.
\end{array}
\end{equation}

For $k\ge k_\circ$, $j=1,\dots,N$, and $t\in[t',t'']$ denote by
$W^j_k(t_1,t_2;t)$ the connected component of the open set $\{x\in
V(t_2)\setminus \overline V(t_1):\Phi_k(x)>-t\}$ such that
$\partial W^j_k(t_1,t_2;t)\supset A_j(t_1)$ (see Fig.1) and put
$$W_k(t_1,t_2;t)=\bigcup\limits_{j=1}^N W^j_k(t_1,t_2;t),\qquad
S_k(t_1,t_2;t)=(\partial W_k(t_1,t_2;t))\cap V(t_2)
\setminus\overline V(t_1).$$ Clearly, $\Phi_k\equiv -t$ on
$S_k(t_1,t_2;t)$. By construction (see Fig.1),
\begin{equation}\label{boundary3}\partial W_k(t_1,t_2;t)=S_k(t_1,t_2;t)\cup A_{1}(t_1)\cup\dots\cup
A_N(t_1).\end{equation} (Note that $W_k(t_1,t_2;t))$ and
$S_k(t_1,t_2;t)$ are well defined for all $t\in[t',t'']$ and $k\ge
k_\circ=k_\circ(t_1,t_2,t',t'')$.)

Since by (E--NS) each $\Phi_k$ belongs to $C^\infty(\Omega_{bk})$,
by the classical Morse-Sard theorem we have that for almost all
$t\in[t',t'']$ the level set $S_k(t_1,t_2;t)$ consists of finitely
many $C^\infty$-cycles and $\Phi_k$ is differentiable (in
classical sense) at every point~$x\in S_k(t_1,t_2;t)$ with
$\nabla\Phi_k(x)\ne0$. The values $t\in[t',t'']$ having the above
property will be called $k$-{\it regular}.

By construction, for every $k$-regular value~$t\in[t',t'']$ the
set $S_{k}(t',t'';t)$ is a finite union of smooth cycles, and
\begin{equation}\label{lac-2-ax}
\int\limits_{S_k(t_1,t_2;t)}\nabla\Phi_k\cdot{\bf
n}\,ds=-\int\limits_{S_k(t_1,t_2;t)}|\nabla\Phi_k|\,ds<0,
\end{equation}
where $\n$ is the unit outward normal vector to $\partial
W_k(t_1,t_2;t)$.

The last inequality leads us to the~main result of this
subsection.

\begin{Lemma}
\label{lem_Leray_fc} {\sl Assume that  $\Omega\subset\R^2$ is
a~domain of type \eqref{Omega} with
$C^2$-smooth boundary $\partial\Omega$,  and~${\bf a}\in
W^{1/2,2}(\partial\Omega)$ satisfies zero total flux
condition~(\ref{dom5}). Then assumptions (E-NS) and
\eqref{as-prev1} lead to a contradiction.}
\end{Lemma}

\proof Fix $t_1,t_2,t',t''\in\Ti$ with $t_1<t'<t''<t_2$. Below we
always assume that $k\ge k_\circ(t_1,t_2,t',t'')$ (see
(\ref{boundary0})--(\ref{boundary2})\,), in particular, the set
$S_k(t_1,t_2;t)$ is well defined for all $t\in[t',t'']$.

The main idea of the proof of Lemma~\ref{lem_Leray_fc} is quite
simple: we will integrate the equation
\begin{equation}\label{cle_lap**}\Delta\Phi_k=\omega_k^2+\frac1{\nu_k}\div(\Phi_k\ue_k)\end{equation} over the suitable
domain $\Omega_{k}(t)$ with $\partial\Omega_{k}(t)\supset
S_k(t_1,t_2;t)$.

We split the construction of the domain $\Omega_{k}(t)$ into
two steps. Namely, for $t\in\Ti\cap[t',t'']$ and sufficiently
large~$k$ denote by $\Omega_{S_k(t_1,t_2;t)}$ the bounded open set
in $\R^2$ such that
$$\partial\Omega_{S_k(t_1,t_2;t)}=  S_k(t_1,t_2;t).$$ Then
put by definition 
\begin{equation}\label{dom-k}
\Omega_{k}(t)=B_k\setminus \Omega_{S_k(t_1,t_2;t)}
\end{equation}
(see Fig.2).
Here $B_k=\{x\in\R^2:|x|<R_k\}$ are the balls where the solutions
$\ue_k\in W^{1,2}(\Omega\cap B_k)$ from (E-NS)-assumptions are
defined.

\begin{center}
\includegraphics[scale=0.4]{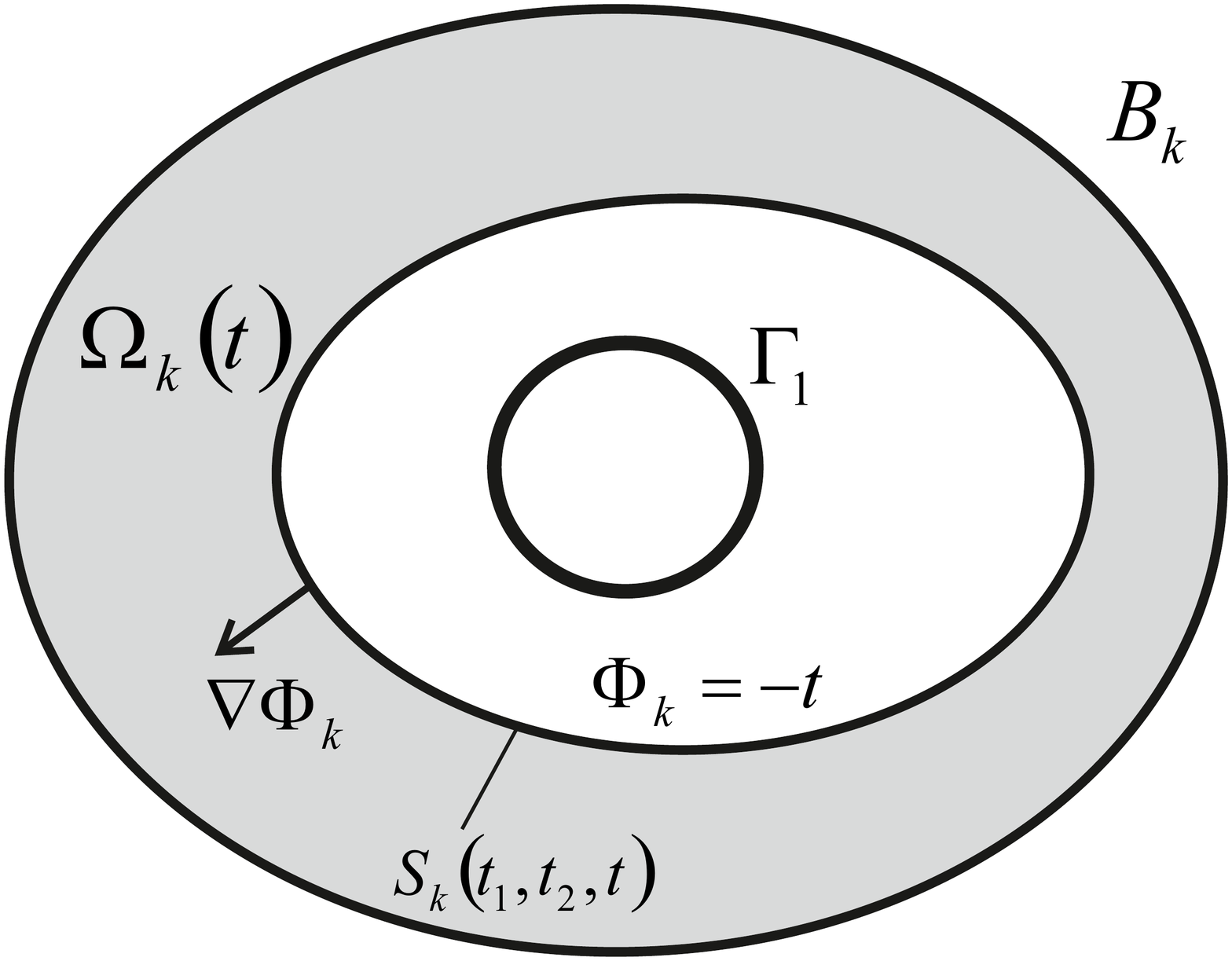}
\end{center}
\begin{center}
Fig. 2. {\sl The domain $\Omega_{k}(t)$ for the case of $N=1$.}
\end{center}

By construction (see Fig.2), \ $\partial\Omega_{k}(t) =
S_k(t_1,t_2;t)\cup C_{R_k }$. Integrating the equation~(\ref{cle_lap**})
over the domain $\Omega_{k}(t)$, we obtain
\begin{equation}\label{cle_lap1}
\begin{array}{lcr}
\int\limits_{S_k(t_1,t_2;t)}\nabla\Phi_k\cdot{\bf n}\,ds+\int\limits_{ C_{R_k
}}\nabla\Phi_k\cdot{\bf n}\,ds=\int\limits_{\Omega_{k}(t)}\omega_k^2\,dx
\\
\\
+\frac{1}{\nu_k}\int\limits_{ S_k(t_1,t_2;t)}\Phi_k\ue_k\cdot{\bf
n}\,ds+\frac{1}{\nu_k}\int\limits_{ C_{R_k}}\Phi_k\ue_k\cdot{\bf n}\,ds.
\end{array}
\end{equation}

By  direct calculations, (\ref{NSk}) implies
\begin{equation}
\label{PP1}\nabla\Phi_k=-\nu_k\nabla^\bot\omega_k+\omega_k\u_k^\bot,
\end{equation}
where, recall, for $\u=(u_1,u_2)$ we denote $\u^\bot=(u_2,-u_1)$
and $\nabla^\bot\omega=(\partial_2\omega,-\partial_1\omega)$.

By the Stokes theorem, for any $C^1$-smooth closed curve
$S\subset\Omega$ and $g\in C^1(\Omega)$ we have
$$\int\limits_S\nabla^\bot g\cdot{\bf n}\,ds=0.$$  So, in particular,
\begin{equation}\label{uesk0}\int\limits_S\nabla\Phi_k\cdot{\bf
n}\,ds=\int\limits_S\omega_k\u_k^\bot\cdot{\bf n}\,ds.
\end{equation}
Since by construction for every $ x\in
C_{R_k}=\{y\in\R^3:|y|=R_k\}$ there holds the equality
\begin{equation}\label{uesk}
\ue_k(x)\equiv \0,
\end{equation}
we see that
\begin{equation}\label{lac0.5}
\int\limits_{C_{R_k}}\nabla\Phi_k\cdot{\bf n}\,ds=0.
\end{equation}

Furthermore, using~(\ref{uesk}) we get
\begin{equation}\label{remop2'}
\frac1{\nu_k}\int\limits_{C_{R_k}}\Phi_k\ue_k\cdot{\bf n}\,ds=0.
\end{equation}

Finally, since $\Phi_k(x)\equiv -t$ for all $x\in S_k(t_1,t_2;t)$,
we obtain
\begin{equation}\label{remop00'}
\int\limits_{S_k(t_1,t_2;t)}\Phi_k\ue_k\cdot{\bf
n}\,ds=-t\int\limits_{S_k(t_1,t_2;t)}\ue_k\cdot{\bf
n}\,ds=t\int\limits_{C_{R_k}}\ue_k\cdot{\bf n}\,ds=0,
\end{equation}
here we have used the identity
$$\int\limits_{\partial\Omega_k(t)}\ue_k\cdot{\bf n}\,ds=
\int\limits_{C_{R_k}}\ue_k\cdot{\bf n}\,ds
+\int\limits_{S_k(t_1,t_2;t)}\ue_k\cdot{\bf n}\,ds=0.$$ In view of
(\ref{lac-2-ax}),  (\ref{cle_lap1}) and (\ref{lac0.5})--(\ref{remop00'})
 we get
\begin{equation}\label{cle_lap2}
\int\limits_{   S_k(t_1,t_2;t)}|\nabla\Phi_k|\,ds =
-\int\limits_{\Omega_{k}(t)}\omega_k^2\,dx,
\end{equation} a contradiction.
The Lemma is proved. \qed

\subsection{The case $\Phi_\infty<\widehat
p_N=\esssup\limits_{x\in\bar\Omega}\Phi(x)$.}\label{Euler-contr-2}

Suppose now that (\ref{as1-axxx}) holds, i.e., the maximum of
$\Phi$ is attained on the boundary component~$\Gamma_N$ and not at
infinity. Then the proof can be reduced to the~case  with
a~bounded domain, which was considered in~\cite{KPRAM}. Let us
describe the essential  details of this reduction.

Without loss of generality we can assume that $\Phi_\infty<0$ and
$\esssup\Phi(x)=\widehat p_N=\Phi(\Gamma_N)=0$. Repeating the
arguments from the first part of
Subsection~\ref{EPcontr-axx}, we construct a~$C^1$~-smooth
cycle $A_N\subset\Omega$ such that $\psi|_{A_N}=\const$, \
$\Phi_\infty<\Phi(A_N)<0$ and $\Gamma_N$ lies in the bounded
connected component of the set~$\R^2\setminus A_N$. Denote this
component by~$\Omega_b$. The cycle $A_N$ separates $ \Gamma_N$
from infinity. Thus, in order
to obtain a~contradiction, it is enough to consider the bounded domain~$\Omega_b\cap\Omega$.

Namely, let $$\Omega_b\cap \Gamma_j=\emptyset,\qquad j=1,\dots
M_1-1,$$
$$\Omega_b\supset \Gamma_j,\qquad j=M_1,\dots,N$$
(the case $M_1=N$ is not excluded). Making a renumeration (if
necessary), we may assume without loss of generality that
$$\Phi( \Gamma_j)<0,\qquad j=M_1,\dots,M_2,$$
$$\Phi( \Gamma_j)=\widehat p_N=0,\qquad j=M_2+1,\dots,N$$
(the case $M_2=M_1-1$, i.e., when $\Phi$ attains maximum value at
every boundary component inside the domain~$\Omega_b$, is not
excluded). Now in order to receive the required contradiction, one
need to~repeat almost word by word the corresponding arguments of
Subsection~2.4.1 in~\cite{KPRAM}. The only modifications are as
follows: now the sets~$A_N$ and $\Gamma_{M_1},\dots, \Gamma_{M_2}$
play the role of the sets~$\Gamma_0,\Gamma_1,\dots,\Gamma_M$
from~\cite[Subsection~2.4.1]{KPRAM}, and  the domain
$\Omega_b\cap\Omega$ in the present case plays the role of the
domain~$\Omega$ from~\cite[Subsection~2.4.1]{KPRAM}, etc.

\subsection{The case $\Phi_\infty=\widehat
p_N=\esssup\limits_{x\in\bar\Omega}\Phi(x)$.}\label{Euler-contr-3}

Consider the last possible case, when the maximum of $\Phi$ is
attained both at infinity and on some boundary component:
\begin{equation}\label{lcas1}
\Phi_\infty=\esssup\limits_{x\in\Omega}\Phi(x)=\widehat
p_N=\max\{\widehat p_1,\dots,\widehat p_N\}
\end{equation}
(recall, that $\widehat p_j=\Phi(\Gamma_j)$\,).

This case is more delicate: we need to combine the arguments of
the previous subsections.

Without loss of generality we may assume that
\begin{equation}\label{lcas2}
0=\Phi_\infty=\esssup\limits_{x\in\Omega}\Phi(x),
\end{equation}
\begin{equation}\label{lcas3}
\pp_j<0,\qquad j=1,\dots, M,
\end{equation}
\begin{equation}\label{lcas4}
\pp_j=0,\qquad j=M+1,\dots, N.
\end{equation}
Note that $1\le M<N$, i.e., the case $\pp_j\equiv0$
for all $j=1,\dots,N$ is impossible. Indeed, from~(\ref{cont_e}) and
(\ref{2.1}$_1$) we have
\begin{equation}\label{lcas5}
-\nu=\sum\limits_{j=1}^N\pp_j\F_j,
\end{equation}
where, recall,
\begin{equation}
\label{SssNS11} \mathscr{F}_j=\into{\partial\Omega_j}\a\cdot\n\,ds.
\end{equation}
Let
$$\delta>\max\{-\pp_j:j=1,\dots,M\}.$$
Using precisely the same arguments as above in
Subsection~\ref{EPcontr-axx}, we construct  a measurable set
$\Ti\subset[0,\delta]$ of full measure (i.e., ${\rm
meas}\,\bigl([0,\delta]\setminus\Ti\bigr)=0$\,) and smooth cycles
$A_j(t)\Subset \Omega$ for all $t\in\Ti$ and every $j=1,\dots,M$
with the following properties:

(i) The set $\Gamma_j=\partial\Omega_j$ lies in a~bounded
connected component of the open set~$\R^2\setminus A_j(t)$;

(ii) $\psi|_{A_j(t)}\equiv\const$, \ $\Phi(A_j(t))=-t$;

(iii) (monotonicity) If $0<t_1<t_2<\delta_{p}$, then $A_j(t_1)$
lies in the unbounded connected component of the
set~$\R^2\setminus A_j(t_2)$ (i.e., the set
$A_j(t_2)\cup\Gamma_j$ lies in the bounded connected component of
the set~$\R^2\setminus A_j(t_1)$);

(iv) $A_j(t)$ is {\it a~regular cycle}, i.e., it is a smooth curve
homeomorphic  to the unit circle and
\begin{equation}
\label{uuu} \mbox{$\Phi_{k}|_{A_j(t)}$ converges to
$\Phi|_{A_j(t)}$ \,uniformly for all $t\in\Ti$.}
\end{equation}

Further, using also the methods of Subsection~\ref{EPcontr-axx},
for any numbers $t_1,t_2,t',t''\\ \in\Ti$ with $t_1<t'<t''<t_2$  and
for all $t\in\Ti\cap(t',t'')$ and $k\ge k_\circ(t_1,t_2,t',t'')$
we construct\footnote{See, e.g.,
 (\ref{boundary0})--(\ref{boundary2}), where now the number
$M$ plays the role of~$N$.} a~domain $\Omega_k(t)$ with
$\partial\Omega_k(t)=C_{R_k}\cup S_k(t_1,t_2;t)$, where
$S_k(t_1,t_2;t)$ is a union of smooth cycles satisfying the
following conditions:
\begin{eqnarray}
\label{uuu1}
\mbox{$S_k(t_1,t_2;t)$ separates $A_j(t_1)$ from $A_j(t_2)$ for all $j\in1,\dots,M$}; \\
\label{uuu2}
\mbox{$\Phi_k\equiv-t$ on $S_k(t_1,t_2;t)$}; \\
\label{uuu3}
\mbox{$\nabla\Phi\ne0$ on $S_k(t_1,t_2;t)$};\\
\label{uuu4}\int\limits_{S_k(t_1,t_2;t)}\nabla\Phi_k\cdot{\bf
n}\,ds=-\int\limits_{S_k(t_1,t_2;t)}|\nabla\Phi_k|\,ds<0,
\end{eqnarray}
where $\n$ is the unit outward normal vector to $\partial
\Omega_k(t)$.

Now we are ready to prove the key estimate.

\begin{Lemma}
\label{ax-lkr11}{\sl For  any $t_1,t_2,t',t''\in\Ti$ with
$t_1<t'<t''<t_2$ there exists $k_*=k_*(t_1,t_2,t',t'')$ such that
for every $k\ge k_*$ and for almost all $t\in[t',t'']$ the
inequality
\begin{equation}\label{mec}
\int\limits_{S_k(t_1,t_2;t)}|\nabla\Phi_k|\,ds<\F t,
\end{equation}
holds with the constant $\F$ independent of  $t,t_1,t_2,t',t''$
and $k$. }
\end{Lemma}

\proof Fix $t_1,t_2,t',t''\in\Ti$ with $t_1<t'<t''<t_2$. Below we
always assume that $k\ge k_\circ=k_\circ(t_1,t_2,t',t'')$, in
particular, the set $S_k(t_1,t_2;t)$ is well defined for all
$t\in[t',t'']\cap\Ti$.

Put $\WQ_k(t)=\Omega\cap\Omega_k(t)$. By construction,
\begin{equation}
\label{uuu7}
\partial \WQ_k(t)=C_{R_k}\cup S_k(t_1,t_2;t)\cup \Gamma_K\cup\dots\cup\Gamma_N,
\end{equation}
where $M<K$. This representation follows from the fact that the
set $S_k(t_1,t_2;t)$ separates the circle $C_{R_k}$ from the
boundary components $\Gamma_j$ with $j=1,\dots,M$. However, a priory
it does not separate $C_{R_k}$  from other boundary components $\Gamma_i$ with
$i>M$. This is the main difference comparing to the situation of
Subsection~\ref{EPcontr-axx}, where the boundary of the
integration domain consists of only two parts: $C_{R_k}\cup
S_k(t_1,t_2;t)$ (see the proof of Lemma~\ref{lem_Leray_fc}\,).

It is easy to see that $K$ in the representation~(\ref{uuu7}) does
not depend on~$k$ for sufficiently large~$k$; see, e.g.,
\cite[Subsection~2.4.1]{KPRAM} for the detailed explanation of
this fact.

Now we have to consider two possible cases:

\medskip

{\sc Case I}. $K=N+1$. It means that no component $\Gamma_j$ is contained in the domain~$ \WQ_k(t)$, i.e.
\begin{equation}
\label{uuu7--}
\partial \WQ_k(t)=C_{R_k}\cup S_k(t_1,t_2;t).
\end{equation}
The contradiction for this case is derived exactly in the same way as in the proof of previous
Lemma~\ref{lem_Leray_fc}.

\

{\sc Case II}. \  $K\le N$. For $h>0$ denote $\Gamma_0=\Gamma_K\cup\dots\cup\Gamma_N$,
$\Gamma_h=\{x\in \Omega:\dist(x,\Gamma_0)=h\}$,
$\Omega_k(t,h)=\{x\in\WQ_k(t):\dist(x,\Gamma_0)>h\}$. Then
\begin{equation}
\label{uuu8}
\partial\Omega_k(t,h)=C_{R_k}\cup S_k(t_1,t_2;t)\cup \Gamma_h
\end{equation}
for any fixed~$t\in\Ti\cap[t',t'']$, for sufficiently
small~$h<\delta(t_1)$ and for sufficiently large~$k\ge k_\circ$.

It was proved in \cite{KPRAM} (see pages~787--788) that for any
fixed~$\e>0$ and for sufficiently large~$k\ge k_\e\ge k_\circ$
there exist a~value $\hh_k<\delta(t_1)$ such that
\begin{equation}\label{iii1}
\biggr|\int\limits_{\Gamma_{\hh_k}}\nabla\Phi_k\cdot{\bf
n}\,ds\biggr|<\e,
\end{equation}
\begin{equation}\label{iii2}
\frac1{\nu_k}\biggr|\int\limits_{\Gamma_{\hh_k}}\Phi_k\ue_k\cdot{\bf
n}\,dS\biggr|<\e.
\end{equation}
It was shown before (see formulas~(\ref{lac0.5})--(\ref{remop2'})\,) that
\begin{equation}\label{iii3}
\int\limits_{C_{R_k}}\nabla\Phi_k\cdot{\bf n}\,ds=0,
\end{equation}
\begin{equation}\label{iii4}
\int\limits_{C_{R_k}}\Phi_k\ue_k\cdot{\bf n}\,ds=0.
\end{equation}
Denote $\Omega_{0k}(t):=\Omega_{k}(t,\hh_k)$. Then
$$\partial\Omega_{0k}(t)=C_{R_k}\cup S_k(t_1,t_2;t)\cup \Gamma_{\hh_k}.$$

Integrating the equation (\ref{cle_lap**}) over the domain $\Omega_{0k}(t)$ and using
(\ref{iii3})--(\ref{iii4}), we get
\begin{equation}\label{lc-cle_lap1}
\begin{array}{lcr}
\int\limits_{S_k(t_1,t_2;t)}\nabla\Phi_k\cdot{\bf
n}\,ds+\int\limits_{\Gamma_{\hh_k}}\nabla\Phi_k\cdot{\bf
n}\,ds=\int\limits_{\Omega_{k}(t)}\omega_k^2\,dx
\\
\\
+\frac{1}{\nu_k}\int\limits_{ S_k(t_1,t_2;t)}\Phi_k\ue_k\cdot{\bf
n}\,ds+\frac{1}{\nu_k}\int\limits_{\Gamma_{\hh_k}}\Phi_k\ue_k\cdot{\bf
n}\,ds.
\end{array}
\end{equation}
Using (\ref{uuu4}), (\ref{iii1})--(\ref{iii2}), we obtain the
estimate
\begin{equation}\label{iii8}
\int\limits_{
S_k(t_1,t_2;t)}|\nabla\Phi_k|\,ds<2\e-\frac{1}{\nu_k}\int\limits_{
S_k(t_1,t_2;t)}\Phi_k\ue_k\cdot{\bf n}\,ds.
\end{equation}

Finally, since $\Phi_k(x)\equiv -t$ for all $x\in S_k(t_1,t_2;t)$,
we derive
\begin{equation}\label{iii9}
\begin{array}{lcr}
\int\limits_{S_k(t_1,t_2;t)}\Phi_k\ue_k\cdot{\bf
n}\,ds=-t\int\limits_{S_k(t_1,t_2;t)}\ue_k\cdot{\bf
n}\,ds=t\int\limits_{\Gamma_0}\ue_k\cdot{\bf
n}\,ds\\
\\
=t\,\nu_k\int\limits_{\Gamma_0}{\bf a}\cdot{\bf
n}\,ds=t\nu_k\F_\circ,
\end{array}
\end{equation}
here $\F_\circ=\frac1\nu\sum\limits_{j=K}^N \F_j$ and we  have used the
identities (\ref{NSk}$_3$), (\ref{uuu7}) and
$$
0=\int\limits_{\partial\widetilde\Omega_k(t)}\ue_k\cdot{\bf n}ds=
\int\limits_{C_{R_k}\cup S_k(t_1,t_2;t)\cup
\Gamma_0}\ue_k\cdot{\bf
n}ds=\int\limits_{S_k(t_1,t_2;t)}\ue_k\cdot{\bf
n}ds+\int\limits_{\Gamma_0}\ue_k\cdot{\bf n}ds.
$$
 Since the parameter
$\e>0$ could be chosen to be arbitrary small, from
 (\ref{iii8})--(\ref{iii9})
it follows the inequality
\begin{equation}\label{iii10}
\int\limits_{   S_k(t_1,t_2;t)}|\nabla\Phi_k|\,ds \le
\bigl(|\F_\circ|+1\bigr)t
\end{equation} for sufficiently large~$k$.
The Lemma is proved. \qed
\\

Now we  apply the argument from \cite[proof of
Lemma~3.9]{KPRAM} and receive the required contradiction using
the~Coarea formula.

\begin{Lemma}
\label{lem_Leray_fc-00} {\sl Assume that  $\Omega\subset\R^2$ is a
bounded domain of type \eqref{Omega} with $C^2$-smooth boundary
$\partial\Omega$,  and ${\bf a}\in W^{1/2,2}(\partial\Omega)$
satisfies  zero total flux condition~(\ref{dom5}). Then
assumptions (E-NS) and \eqref{as-prev-id-ax} lead to a
contradiction.}
\end{Lemma}

\proof Take a number $t_0\in\Ti$ such that $t_i:=2^{-i}t_0\in\Ti$
for all~$i\in\N$. Let $R_0$  be a~sufficiently large radius such that
$B_{\frac12 R_0}\supset\partial\Omega$. Denote $S_{ik}(t):=
B_{R_0}\cap S_{k}(t_{i+1},t_i,\frac58t_i,\frac78t_i)$ (it is well
defined for almost all $t\in[\frac58t_i,\frac78t_i]$ and for $k\ge
k_*\ge k_\circ$, see paragraph before Lemma~\ref{ax-lkr11}\,) and
put
$$
E_i=\bigcup\limits_{t\in [\frac58t_i,\frac78t_i]\cap\Ti} S_{ik}(t).
$$
 By the Coarea formula (see, e.g.,~\cite{Maly}), for any integrable
 function $g:E_i\to\R$    the equality
\begin{equation}\label{Coarea_Phi}\int\limits_{E_i}g|\nabla\Phi_k|\,dx=
\int\limits_{\frac58t_i}^{\frac78t_i}\int_{
S_{ik}(t)}g(x)\,ds\,dt
\end{equation}
holds.  In particular, taking $g=|\nabla\Phi_k|$ and
using~(\ref{mec}), we obtain
\begin{equation}\label{Coarea_Phi2}\int\limits_{E_i}|\nabla\Phi_k|^2\,dx=
\int\limits_{\frac58t_i}^{\frac78t_i}\int_{
S_{ik}(t)}|\nabla\Phi_k|(x)\,ds\,dt \le
\int\limits_{\frac58t_i}^{\frac78t_i}\F t\,dt= \F't_i^2
\end{equation}
where $\F'=\frac3{16}\F$ is independent of $i$. Now, taking  $g=1$
in (\ref{Coarea_Phi})  and using the H\"older inequality we have
\begin{equation}\label{Coarea_Phi3}
\begin{array}{lcr}
\displaystyle \int\limits_{\frac58t_i}^{\frac78t_i}\mathfrak{H}^1\bigl(
S_{ik}(t)\bigr)\,dt= \int\limits_{E_i}|\nabla\Phi_k|\,dx
\\
\displaystyle \le
\biggl(\int\limits_{E_i}|\nabla\Phi_k|^2\,dx\biggr)^{\frac12}
\bigl(\meas (E_i)\bigr)^{\frac12}\le\sqrt{\F'}t_i\bigl(\meas
(E_i)\bigr)^{\frac12}.
\end{array}
\end{equation}
 By construction, for almost all $t\in [\frac58t_i,\frac78t_i]$ the set $S_{ik}(t)$ is a finite union of smooth lines and
$S_{ik}(t)$ separates $A_j(t_{i+1})$ from $A_j(t_i)$  \, in
$B_{R_0}$ for $j=1,\dots,M$. Thus, each set $S_{ik}(t)$ separates
$\Gamma_j$ from $\Gamma_N$. In particular,
$\Ha^1(S_{ik}(t))\\
\ge\min\bigl(\diam(\Gamma_j),\diam(\Gamma_N)\bigr)$.
Hence, the left integral in (\ref{Coarea_Phi3}) is greater than
$Ct_i$, where $C>0$ does not depend on $i$. On the other hand, the
sets $E_i$ are  pairwise disjoint and, therefore, $\meas(E_i)\to0$ as
$i\to\infty$. The obtained contradiction finishes the proof of
Lemma~\ref{lem_Leray_fc-00}.
$\qed$\\

 We can summarize the results of Subsections~\ref{EPcontr-axx}--\ref{Euler-contr-3} in the following statement.

\begin{Lemma}
\label{lem_Leray_fc_3} {\sl Assume that $\Omega\subset\R^2$ is an
exterior  plane domain of type~\eqref{Omega} with $C^2$-smooth
boundary $\partial\Omega$ and ${\bf
a}\in W^{1/2,2}(\partial\Omega)$ satisfies  zero total flux
condition~(\ref{dom5}). Let (E-NS) be fulfilled. Then every
possible assumption (\ref{as-prev1}), (\ref{as1-axxx}) and
(\ref{as-prev-id-ax}) lead to a contradiction. }
\end{Lemma}

 {\bf Proof of Theorem \ref{T2}.} Let the hypotheses
of Theorem \ref{T2} be satisfied. Suppose that its assertion
fails. Then, by Lemma~\ref{lem_Leray_symm}, there exist $ \ve, p$
and a sequence $(\ue_k,p_k)$ satisfying (E-NS), and by
Lemmas~\ref{lem_Leray_fc_3} these assumptions lead to a
contradiction.    \qed
\\

 {\em {\bf Acknowledgment.}} M.~Korobkov was partially supported by
the Ministry of Education and Science of the Russian Federation
(grant 14.Z50.31.0037). The main part of the paper was written
during a visit of M.~Korobkov to the University of Campania "Luigi
Vanvitelli" in~2017, and he is very thankful for the hospitality.

The research of K. Pileckas was funded by a grant No. S-MIP-17-68 from the Research Council of Lithuania.

{\small 
\end{document}